\documentclass[11pt,reqno]{amsart}

\usepackage[a4paper,margin=1.15in]{geometry}
\usepackage{amsmath,amssymb,amsfonts,amsthm,mathtools}
\usepackage{mathrsfs}
\usepackage{enumitem}
\usepackage{hyperref}
\usepackage{xcolor}
\usepackage{booktabs}
\usepackage{aliascnt}

\hypersetup{
  colorlinks=true,
  linkcolor=blue!55!black,
  citecolor=blue!55!black,
  urlcolor=blue!55!black,
  pdftitle={Quantitative Driver-Only Capacity Cores and Simultaneous Scalar Ito Flows},
  pdfauthor={Guangqian Zhao}
}

\theoremstyle{plain}
\newtheorem{theorem}{Theorem}[section]
\newaliascnt{proposition}{theorem}
\newtheorem{proposition}[proposition]{Proposition}
\aliascntresetthe{proposition}
\newaliascnt{lemma}{theorem}
\newtheorem{lemma}[lemma]{Lemma}
\aliascntresetthe{lemma}
\newaliascnt{corollary}{theorem}
\newtheorem{corollary}[corollary]{Corollary}
\aliascntresetthe{corollary}

\theoremstyle{definition}
\newaliascnt{definition}{theorem}
\newtheorem{definition}[definition]{Definition}
\aliascntresetthe{definition}
\newaliascnt{assumption}{theorem}

\aliascntresetthe{assumption}
\newaliascnt{remark}{theorem}
\newtheorem{remark}[remark]{Remark}
\aliascntresetthe{remark}
\newaliascnt{example}{theorem}
\newtheorem{example}[example]{Example}
\aliascntresetthe{example}
\newaliascnt{problem}{theorem}
\newtheorem{problem}[problem]{Problem}
\aliascntresetthe{problem}

\usepackage[nameinlink,capitalize]{cleveref}
\crefname{theorem}{Theorem}{Theorems}
\crefname{proposition}{Proposition}{Propositions}
\crefname{lemma}{Lemma}{Lemmas}
\crefname{corollary}{Corollary}{Corollaries}
\crefname{definition}{Definition}{Definitions}
\crefname{assumption}{Assumption}{Assumptions}
\crefname{remark}{Remark}{Remarks}
\crefname{example}{Example}{Examples}
\crefname{problem}{Problem}{Problems}
\Crefname{theorem}{Theorem}{Theorems}
\Crefname{proposition}{Proposition}{Propositions}
\Crefname{lemma}{Lemma}{Lemmas}
\Crefname{corollary}{Corollary}{Corollaries}
\Crefname{definition}{Definition}{Definitions}
\Crefname{assumption}{Assumption}{Assumptions}
\Crefname{remark}{Remark}{Remarks}
\Crefname{example}{Example}{Examples}
\Crefname{problem}{Problem}{Problems}

\newcommand{\R}{\mathbb R}
\newcommand{\N}{\mathbb N}

\newcommand{\E}{\mathbb E}
\newcommand{\cB}{\mathcal B}
\newcommand{\cC}{\mathcal C}
\newcommand{\cD}{\mathcal D}
\newcommand{\cF}{\mathcal F}
\newcommand{\cG}{\mathcal G}
\newcommand{\cI}{\mathcal I}
\newcommand{\cK}{\mathcal K}
\newcommand{\cP}{\mathcal P}
\newcommand{\cR}{\mathcal R}

\newcommand{\eps}{\varepsilon}
\newcommand{\one}{\mathbf 1}
\newcommand{\dd}{\,d}
\newcommand{\qv}{\mathcal Q}

\newcommand{\Var}{\operatorname{Var}}

\newcommand{\Graph}{\operatorname{Graph}}
\newcommand{\norm}[1]{\left\lVert #1\right\rVert}
\newcommand{\abs}[1]{\left\lvert #1\right\rvert}

\title[Quantitative driver-only capacity cores]
{Quantitative Driver-Only Capacity Cores and Simultaneous Scalar It\^o Flows}
\author{Guangqian Zhao}
\address{School of Mathematical Sciences, University of Science and Technology of China, Hefei, Anhui 230026, China}
\email{zhaoguangqian@mail.ustc.edu.cn}


\begin{document}
\raggedbottom

\begin{abstract}
We construct quantitative compact raw-driver cores for the maximal nondominated class $M_\Lambda$ of laws under which the coordinate process is a scalar continuous local martingale satisfying $d[X]_t\leq\Lambda\,dt$. Fix $0<\eta<1/2$ and $0<\gamma<\min\{1/2,2\eta\}$. Explicit compact sets $K_{A,B}$ have sub-Gaussian tails in $A$ and subexponential tails in $B$, while continuous dyadic realized variation converges uniformly on each core at rate $O(B\Lambda T2^{-\gamma n})$. The intrinsic dyadic quadratic variation is Holder continuous in the raw uniform topology on the cores. The parabolic family $K_R=K_{R,R^2}$ has complement of capacity at most $Ce^{-cR^2}$ and is stable under stopping and controlled deterministic concatenation.

For every compact family of sufficiently smooth autonomous coefficients with uniformly positive diffusion, a deterministic Lamperti--Follmer equation constructs Borel causal solution and integral fields jointly in the coefficient, starting time, and initial value. The fields satisfy the solution and integral cocycles. The solution field forms an orientation-preserving $C^1$ flow, and both fields have quantitative raw-driver Holder moduli on each $K_R$. Under every law in $M_\Lambda$, each fixed section agrees with the corresponding classical Ito solution and stochastic integral.

We also place the Bichteler--Karandikar construction on a fixed Borel pair-domain, with projective causality and full-sequence semimartingale identification.
\end{abstract}

\maketitle

\section{Introduction}\label{sec:introduction}

\subsection{The pathwise question}

Let \(\Omega\) be a canonical path space and let \(\cP\) be a family of
probability laws on \(\Omega\), not assumed to be dominated.  A stochastic
equation may have the same formal expression under every \(P\in\cP\), and may
be well posed \(P\)-almost surely for each \(P\).  This does not by itself
produce a single pathwise equation or a single solution map.  The quantifiers
\[
  \forall P\in\cP\;\exists N_P\quad P(N_P)=0
\]
do not identify a common process or a fixed failure event \(N\) such that the
same pathwise property holds on \(N^c\).  In particular, the union of the
measure-dependent exceptional sets need not be polar.

One can instead begin with a Borel functional defined by a common
deterministic approximation sequence.  Its convergence event is a fixed set
\(G\subset\Omega\).  If
\[
  P(G)=1\qquad\text{for every }P\in\cP,
\]
then \(G^c\) is polar, while the value of the functional on \(G\) is
determined by the path.

We implement this principle for an It\^o functional based on endogenous
partitions.  The construction separates the deterministic definition and
causal structure from the probabilistic verification of its domain.
Passing from measurewise solutions to a universal solution map requires an
aggregation or selection argument.

\subsection{Endogenous partitions}

For a c\`adl\`ag path \(h\) and \(\eps_n=2^{-n}\), set
\[
  \tau_0^n(h)=0,\qquad
  \tau_{k+1}^n(h)
  =
  \inf\bigl\{t>\tau_k^n(h):
       \abs{h(t)-h(\tau_k^n(h))}>\eps_n\bigr\}\wedge T.
\]
For a continuous path \(x\), consider
\begin{equation}\label{eq:intro-sum}
  I_t^n(h,x)
  =
  \sum_{k\ge0}h(\tau_k^n(h))
  \bigl(
    x(t\wedge\tau_{k+1}^n(h))
    -
    x(t\wedge\tau_k^n(h))
  \bigr).
\end{equation}
The partition is generated by the integrand and follows its oscillations.
This is the pathwise stochastic-integration construction of Bichteler and
Karandikar \cite{Bichteler81,Karandikar95}.  We use it as a canonical
common-domain mechanism and make its measurable, projective, and causal
structure explicit.

\subsection{Main results}

The main results are as follows.

\begin{enumerate}[leftmargin=2.2em,label=(\roman*)]
  \item We place the Bichteler--Karandikar sums on fixed Borel pair-domains.
  Their partial limits
  are causal and projectively consistent, and the full sequence converges
  almost surely under every continuous-semimartingale law on the joint
  canonical space.  Coordinate integrands then give a Borel
  quadratic-covariation functional and a law-free solution relation.  Standard
  selection results give conditions under which such a relation yields a Borel
  causal solution map.

  \item The main quantitative step uses continuous dyadic realized variation.
  The \(L^r\)-norm of its uniform-in-time error is \(O(r)\), uniformly over
  the maximal bounded-volatility martingale class.  Combining the resulting
  subexponential defect with the sub-Gaussian H\"older
  norm of the driver produces explicit two-scale compact raw-driver cores.
  The parabolic diagonal of this family has Gaussian capacity tails.  On every core the intrinsic
  quadratic variation has an explicit raw-uniform H\"older modulus, and the
  cores are stable under stopping and linearly enlarged concatenation.

  \item A deterministic Lamperti--F\"ollmer equation over the first-level path
  and its intrinsic dyadic quadratic variation constructs one solution field
  and one integral field, simultaneously in the coefficient, starting time,
  and initial value.  Both fields have a quantitative raw-driver modulus on
  the same cores, while the solution field is an orientation-preserving
  \(C^1\) flow.
\end{enumerate}

\begin{theorem}[Main theorem]
\label{thm:main-intro}
Let \(\mathfrak M_\Lambda\) be the maximal class of all laws on
\(C_0([0,T])\) under which the coordinate process is a continuous local
martingale with \(\dd[X]_t\le\Lambda\,\dd t\), and set
\(c_\Lambda(A)=\sup_{P\in\mathfrak M_\Lambda}P(A)\).  Fix
\[
  0<\eta<\frac12,
  \qquad
  0<\gamma<\min\left\{\frac12,2\eta\right\},
  \qquad
  \alpha=\frac{\gamma}{1-\eta+\gamma}.
\]
Then the following statements hold.
\begin{enumerate}[label=(\roman*),leftmargin=2.2em]
  \item For every \(A,B\ge1\) there is an explicit compact raw-driver set
  \(\cK_{A,B}\subset\cG_{\mathrm{dyad}}\), and there are positive constants
  depending only on \((\eta,\gamma)\), such that
  \begin{align*}
    c_\Lambda(\cK_{A,B}^c)
      &\le C_\eta e^{-c_\eta A^2}+C_\gamma e^{-c_\gamma B},\\
    \sup_{x\in\cK_{A,B}}\norm{Q^n(x)-Q(x)}_\infty
      &\le C_\gamma B\Lambda T2^{-\gamma n},\\
    \norm{Q(x)-Q(\widetilde x)}_\infty
      &\le C_{\eta,\gamma}A^\alpha B^{1-\alpha}
       (\Lambda T)^{1-\alpha/2}
       \norm{x-\widetilde x}_\infty^\alpha,
       \qquad x,\widetilde x\in\cK_{A,B}.
  \end{align*}
  Stopping preserves \(\cK_{A,B}\), and deterministic concatenation maps
  \(\cK_{A,B}\times\cK_{A,B}\) into
  \(\cK_{2^{1-\eta}A,\,3B+\kappa_\eta A^2}\) for a constant
  \(\kappa_\eta\).  On the parabolic diagonal
  \(\cK_R:=\cK_{R,R^2}\),
  \[
    c_\Lambda(\cK_R^c)\le Ce^{-cR^2},
    \qquad
    x,\widetilde x\in\cK_R
    \ \Longrightarrow\ 
    x\otimes_s\widetilde x\in\cK_{C_*R}.
  \]
  Moreover, for every \(P\in\mathfrak M_\Lambda\),
  \[
    P(\cG_{\mathrm{dyad}})=1,
    \qquad
    Q(X)=[X]^P\quad P\text{-almost surely}.
  \]

  \item For every compact coefficient class
  \(\Theta\subset C_b^1(\R)\times C_b^2(\R)\) satisfying, for some
  \(L<\infty\) and
  \(0<\underline\sigma\le\overline\sigma<\infty\),
  \[
    \sup_{\theta\in\Theta}
    \left(
      \norm{b_\theta}_{C_b^1}+\norm{\sigma_\theta}_{C_b^2}
    \right)\le L,
    \qquad
    \underline\sigma\le\sigma_\theta(y)\le\overline\sigma
    \quad(\theta\in\Theta,\ y\in\R),
  \]
  one Borel causal solution field \(S\) and one Borel causal integral field
  \(J\) are defined on the
  fixed Borel domain \(\cG_{\mathrm{dyad}}\), whose complement is
  \(c_\Lambda\)-polar, simultaneously in the coefficient, starting time,
  and initial value.  They satisfy the solution and integral cocycles, the
  solution field is an orientation-preserving \(C^1\) flow, and on every
  \(\cK_R\) both fields are H\"older in the raw driver with exponent
  \[
    \beta_{\mathrm{flow}}
    =\eta\alpha
    =\frac{\eta\gamma}{1-\eta+\gamma}.
  \]
  Under each \(P\in\mathfrak M_\Lambda\), every fixed parameter section is
  the corresponding classical It\^o solution and integral version.
\end{enumerate}
\end{theorem}

Part~(i) is proved in
\Cref{thm:uniform-realized-qv,thm:driver-core,cor:parabolic-core};
part~(ii) is \Cref{thm:simultaneous-flow}.  The core estimates are
deterministic; their capacity bounds and stochastic identifications are
probabilistic.

\subsection{Position relative to existing theories}

Pathwise algorithms for stochastic integration were developed by Bichteler
\cite{Bichteler81} and Karandikar \cite{Karandikar95}.  For nondominated
families, Nutz \cite{Nutz12} constructed simultaneous stochastic integrals
under a common predictable domination of the characteristics, while Soner,
Touzi and Zhang \cite{STZ11} developed an aggregation framework and analyzed
its obstructions.  Related measurable constructions appear in Nutz and van
Handel \cite{NutzVanHandel13}.  Universal solution functionals were studied
by Kallenberg \cite{Kallenberg96}, Kurtz \cite{Kurtz14}, Oberhauser
\cite{Oberhauser16}, and Przyby\l owicz, Schwarz, Steinicke and Sz\"olgyenyi
\cite{PrzybylowiczEtAl24}.  These works provide the aggregation and
pair-domain framework used in \Cref{sec:fixedsets,sec:integral}.

Model-free integration and differential equations have also been developed
through outer-measure and superhedging methods
\cite{LochowskiPerkowskiPromel18,BartlKupperNeufeld19,
LochowskiPerkowskiPromel22,GalaneEtAl23}.  Quasi-continuity under weakly
compact families is treated systematically by Denis, Hu and Peng
\cite{DenisHuPeng11}.  The quantitative object considered here is an
explicit family of compact sets in the raw uniform topology, derived from the
single characteristic bound \(\dd[X]\le\Lambda\,\dd t\).  On these sets the
intrinsic quadratic variation has a raw-path modulus, the approximation rate
is uniform, and stopping and concatenation preserve the core up to controlled
enlargement.  The same cores support one field indexed by the coefficient,
starting time, and initial value.

F\"ollmer's deterministic It\^o formula \cite{Follmer81} supplies the
pathwise change-of-variables step.  Integral equations in It\^o--F\"ollmer
calculus were studied by Hirai \cite{Hirai19}; scalar pathwise transformations
for stochastic integral equations were developed by Karatzas and Ruf
\cite{KaratzasRuf16}.  In the \(G\)-Brownian framework, Luo and Wang
\cite{LuoWang14} represented solutions of one-dimensional \(G\)-SDEs
through ordinary differential equations parametrized by the sample path, with
quadratic variation entering as a finite-variation input.  The theorem below
is instead formulated for the maximal bounded-volatility martingale class and
couples the scalar transformation with explicit compact capacity cores,
raw-path moduli, stopping and concatenation stability, and a single Borel
causal field jointly indexed by the coefficient, starting time, and initial
value.  For classical stochastic-flow theory, see Kunita \cite{Kunita90}.

Rough path theory \cite{Lyons98,FrizHairer20} obtains deterministic
continuity after enhancing the driver.  Quasi-sure rough lifts for
\(G\)-Brownian motion and their relation to \(G\)-SDEs were developed in
\cite{GQY14,PengZhang17}; homeomorphic \(G\)-SDE flows were established by
Gao \cite{Gao09}.  Das, Kwossek and Pr\"omel \cite{DasKwossekPromel25}
related general F\"ollmer-type Riemann sums to non-geometric rough
paths.  The endogenous partitions in \eqref{eq:intro-sum} instead depend on
the integrand, whereas the driver-only theorem uses intrinsic dyadic
quadratic variation.  See \Cref{sec:geometry}.

\subsection{Organization}

\Cref{sec:fixedsets} separates fixed pathwise domains from measurewise
patching and proves the common-approximation principle.
\Cref{sec:integral} constructs the endogenous-partition It\^o functional and
identifies it under arbitrary continuous-semimartingale laws.
\Cref{sec:equations} constructs quadratic covariation and the law-free
solution relation, then proves the aggregation and causality theorems.
\Cref{sec:driver-core} proves the driver-only compact-core and universal flow
theorem.  \Cref{sec:geometry} compares the constructions with enhanced-path
calculi and records the remaining general capacity-core problem.

\section{Fixed pathwise domains and aggregation}
\label{sec:fixedsets}

\subsection{Capacity and fixed exceptional sets}

Let \(E\) be a standard Borel space and let \(\cP\) be a nonempty family of
Borel probability measures on \(E\).  For an arbitrary set \(A\subset E\), set
\begin{equation}\label{eq:capacity}
  c_{\cP}(A):=\sup_{P\in\cP}P^*(A),
\end{equation}
where \(P^*\) is outer probability.  A set \(N\) is
\(\cP\)-\emph{polar} if \(c_{\cP}(N)=0\).  A property holds
\(\cP\)-quasi surely if it holds outside a polar set.  This is the standard
capacity language of nondominated analysis; see
\cite{DenisHuPeng11,STZ11}.

The following elementary facts isolate the role of a common event.

\begin{proposition}[Common full-measure event]\label{prop:fixed-event}
Let \(G\subset E\) be universally measurable.  If
\[
  P(G)=1\qquad\text{for every }P\in\cP,
\]
then \(G^c\) is \(\cP\)-polar.
\end{proposition}

\begin{proof}
For every \(P\in\cP\), universal measurability gives
\(P^*(G^c)=P(G^c)=0\).  Taking the supremum over \(P\) proves the claim.
\end{proof}

\begin{remark}\label{rem:capacity-one}
The statement \(c_{\cP}(G)=1\) is not the appropriate common full-measure
condition: it only asserts that the supremum of \(P(G)\) is one.  The common
condition is
\[
  c_{\cP}(G^c)=0,
  \qquad\text{equivalently}\qquad
  P(G)=1\quad\text{for every }P\in\cP .
\]
\end{remark}

\subsection{An aggregation obstruction}

The following static example exhibits the aggregation obstruction.

\begin{example}[Pairwise compatibility without an aggregator]
\label{ex:no-aggregate}
Let \(E=[0,1]^2\), let \(\lambda\) be Lebesgue measure on \([0,1]\), and for
\(a\in[0,1]\) define
\[
  P_a=\frac12\bigl(\delta_a\otimes\lambda
                    +\lambda\otimes\delta_a\bigr),
  \qquad
  Z^a\equiv a .
\]
For \(a\ne b\), the supports
\[
  S_a=\bigl(\{a\}\times[0,1]\bigr)
      \cup\bigl([0,1]\times\{a\}\bigr)
\]
and \(S_b\) meet in only finitely many points.  These points are null under
both \(P_a\) and \(P_b\).  More explicitly,
\[
  G_a^{a,b}=S_a\setminus S_b,
  \qquad
  G_b^{a,b}=S_b\setminus S_a
\]
are respectively \(P_a\)-full and \(P_b\)-full, and their intersection is
empty.  Hence each pair \(Z^a,Z^b\) can be made compatible on pair-dependent
full-measure subsets.

Here universally measurable has its standard meaning with respect to
\[
  \mathcal U(E)
  =
  \bigcap_{\mu\in\mathfrak P(E)}\cB(E)^\mu,
\]
the intersection of the completions under all Borel probability measures.
There is no universally measurable
\(Z:E\to[0,1]\) such that
\[
  Z=Z^a\qquad P_a\text{-almost surely for every }a\in[0,1].
\]
Indeed, for every \(a\),
\begin{align*}
  1=P_a(Z=a)
  ={}&\frac12\lambda\{y:Z(a,y)=a\}\\
     &+\frac12\lambda\{x:Z(x,a)=a\}.
\end{align*}
Both terms on the right are at most one, so each of the two displayed
Lebesgue measures equals one.  Thus such a \(Z\) would satisfy
\[
  Z(a,y)=a\quad\text{for \(\lambda\)-almost every \(y\)}
\]
and
\[
  Z(x,a)=a\quad\text{for \(\lambda\)-almost every \(x\)}
\]
for every \(a\).  Define
\[
  A=\{(x,y):Z(x,y)=x\},
  \qquad
  B=\{(x,y):Z(x,y)=y\}.
\]
Universal measurability makes \(A\) and \(B\)
\(\lambda^2\)-measurable.  Fubini gives
\(\lambda^2(A)=\lambda^2(B)=1\), whereas
\[
  A\cap B\subset\{(x,y):x=y\},
\]
and the diagonal has zero \(\lambda^2\)-measure, a contradiction.
\end{example}

\begin{remark}
\Cref{ex:no-aggregate} is a static analogue of the aggregation obstruction in
\cite[Example 3.3]{STZ11}.  It shows that a statement of the
form
\[
  \forall P,Q\ \exists\text{ pair-dependent compatible full sets}
\]
cannot be exchanged for
\[
  \exists\text{ a measurable object compatible with every }P .
\]
\end{remark}

\subsection{Common deterministic approximations}

The fixed-event mechanism becomes useful when the event is the convergence
set of one approximation sequence.

\begin{theorem}[Common-approximation principle]\label{thm:common-approx}
Let \(E\) be a standard Borel space, let \((F,d_F)\) be a Polish space, and let
\(\Phi_n:E\to F\) be Borel maps.  Define
\begin{equation}\label{eq:common-convergence-set}
  G_\Phi
  :=
  \{e\in E:(\Phi_n(e))_{n\ge1}\text{ converges in }F\}.
\end{equation}
Then \(G_\Phi\) is Borel and the limit
\[
  \Phi(e):=\lim_{n\to\infty}\Phi_n(e),\qquad e\in G_\Phi,
\]
is Borel as a map from \(G_\Phi\) to \(F\).

If, for every \(P\in\cP\), there is an \(F\)-valued random variable
\(Y^P\) such that
\[
  \Phi_n\longrightarrow Y^P
  \qquad P\text{-almost surely},
\]
then \(G_\Phi^c\) is \(\cP\)-polar and \(\Phi=Y^P\),
\(P\)-almost surely, for every \(P\in\cP\).
\end{theorem}

\begin{proof}
Completeness of \(F\) gives the Cauchy representation
\[
  G_\Phi
  =
  \bigcap_{\ell\ge1}\;
  \bigcup_{N\ge1}\;
  \bigcap_{n,r\ge N}
  \bigl\{e:d_F(\Phi_n(e),\Phi_r(e))<\ell^{-1}\bigr\},
\]
which is Borel.  A pointwise limit of Borel maps into a Polish space is Borel
on its convergence domain.  The almost sure convergence assumption implies
\(P(G_\Phi)=1\) and identifies the limit with \(Y^P\) for every
\(P\).  \Cref{prop:fixed-event} completes the proof.
\end{proof}

\begin{proposition}[Capacity-summable convergence]\label{prop:capacity-bc}
Let \(\Phi_n:E\to F\) be Borel, and let \((\delta_n)\) be a positive sequence
with \(\sum_n\delta_n<\infty\).  If
\begin{equation}\label{eq:capacity-summable}
  \sum_{n=1}^{\infty}
  c_{\cP}\bigl(
    d_F(\Phi_{n+1},\Phi_n)>\delta_n
  \bigr)<\infty,
\end{equation}
then \((\Phi_n)\) converges in \(F\) outside one Borel
\(\cP\)-polar set.
\end{proposition}

\begin{proof}
Let
\[
  A_n=\{d_F(\Phi_{n+1},\Phi_n)>\delta_n\}.
\]
Countable subadditivity of the upper capacity gives
\[
  c_{\cP}\Bigl(\bigcup_{n\ge N}A_n\Bigr)
  \le\sum_{n\ge N}c_{\cP}(A_n).
\]
Hence \(c_{\cP}(\limsup_nA_n)=0\).  Outside the limsup, the increments are
eventually bounded by \(\delta_n\), so the sequence is Cauchy and converges by
completeness of \(F\).
\end{proof}

\section{The endogenous-partition It\^o functional}
\label{sec:integral}

\subsection{Canonical spaces and level-crossing times}

Fix \(T>0\) and integers \(d,m\ge1\).  We use the Polish spaces
\[
  \mathsf X_d=C([0,T];\R^d),\qquad
  \mathsf H_{m,d}=D([0,T];\R^{m\times d}),
\]
where \(\mathsf X_d\) carries the uniform topology and
\(\mathsf H_{m,d}\) carries the Skorokhod \(J_1\)-topology.  Set
\[
  \mathsf E_{m,d}:=\mathsf H_{m,d}\times\mathsf X_d.
\]
We write a generic point as \((h,x)\) and equip \(\mathsf E_{m,d}\) with the
raw canonical filtration
\[
  \cF_t^0
  =
  \sigma\{h(s),x(s):0\le s\le t\}.
\]
Its right-continuous modification is denoted by
\[
  \cF_{t+}^0
  :=
  \bigcap_{\substack{u>t\\u\in\mathbb Q}}\cF_u^0,
  \qquad t<T,
  \qquad
  \cF_{T+}^0=\cF_T^0,
\]
where the intersection is taken over rational \(u\le T\).
All matrix norms below are Euclidean norms after a fixed identification with
a finite-dimensional vector space.

For \(n\ge1\), let \(\eps_n=2^{-n}\) and define
\begin{align}
  \tau_0^n(h)&=0, \label{eq:tau-zero}\\
  \tau_{k+1}^n(h)
  &=
  \inf\bigl\{
      t>\tau_k^n(h):
      \norm{h(t)-h(\tau_k^n(h))}>\eps_n
  \bigr\}\wedge T . \label{eq:tau-recursion}
\end{align}
After the first index \(k\) for which \(\tau_k^n(h)=T\), we set
\(\tau_j^n(h)=T\) for every \(j\ge k\).  The strict inequality in
\eqref{eq:tau-recursion} makes this a first entrance into an open set.

\begin{lemma}[Endogenous crossing times]\label{lem:crossing-times}
For every \(n\ge1\) and \(k\ge0\), the map
\[
  h\longmapsto\tau_k^n(h)
\]
is Borel from \(\mathsf H_{m,d}\) to \([0,T]\).  It is a stopping time for
\((\cF_{t+}^0)_{t\le T}\).  Moreover, for each fixed
\(h\in\mathsf H_{m,d}\), only finitely many \(\tau_k^n(h)\) are strictly
smaller than \(T\).
\end{lemma}

\begin{proof}
We argue simultaneously by induction that \(\tau_k^n\) is Borel and is an
\((\cF_{t+}^0)_{t\le T}\)-stopping time.  This is immediate for
\(\tau_0^n=0\).  Suppose it holds for \(\tau_k^n\).  The joint evaluation map
\[
  (h,t)\longmapsto h(t)
\]
is Borel on
\(D([0,T];\R^{m\times d})\times[0,T]\), so the map
\(h\mapsto h(\tau_k^n(h))\) is Borel.  Define
\[
  Z_t^{n,k}(h)
  :=
  \one_{\{\tau_k^n(h)\le t\}}
  \bigl(h(t)-h(\tau_k^n(h))\bigr).
\]
Then \(Z^{n,k}\) is an \((\cF_{t+}^0)\)-adapted c\`adl\`ag process and
\[
  \tau_{k+1}^n
  =
  \inf\{t\in[0,T]:\norm{Z_t^{n,k}}>\eps_n\}\wedge T.
\]
For every \(r<T\), right-continuity and openness of
\(\{z:\norm z>\eps_n\}\) give
\[
  \{\tau_{k+1}^n<r\}
  =
  \bigcup_{q\in\mathbb Q\cap[0,r)}
  \{\norm{Z_q^{n,k}}>\eps_n\}
  \in\cF_{r+}^0.
\]
Right-continuity of the filtration yields
\[
  \{\tau_{k+1}^n\le r\}
  =
  \bigcap_{\substack{\ell\ge1\\r+\ell^{-1}<T}}
  \{\tau_{k+1}^n<r+\ell^{-1}\}
  \in\cF_{r+}^0.
\]
The assertion at \(r=T\) is immediate.  Thus
\(\tau_{k+1}^n\) is a stopping time.  Its rational sublevel sets are Borel,
so it is also a Borel path function.

If \(\tau_{k+1}^n<T\), choose
\(s_j\downarrow\tau_{k+1}^n\) such that
\[
  \norm{h(s_j)-h(\tau_k^n)}>\eps_n.
\]
Right-continuity gives
\[
  \norm{h(\tau_{k+1}^n)-h(\tau_k^n)}\ge\eps_n.
\]
Right-continuity at \(\tau_k^n\) also implies
\(\tau_{k+1}^n>\tau_k^n\) whenever \(\tau_k^n<T\).  If infinitely many
crossings occurred before \(T\), then \(\tau_k^n\uparrow s\le T\).
The case \(s=0\) is excluded by right-continuity at zero.  For \(s>0\), the
existence of \(h(s-)\) would make \((h(\tau_k^n))_k\) Cauchy, contradicting
the preceding separation.
\end{proof}

\begin{remark}\label{rem:hitting-not-continuous}
The maps \(h\mapsto\tau_k^n(h)\) are generally not continuous in the
\(J_1\)-topology, but they are Borel.  The stopped sums satisfy the
prefix-causality property used below.
\end{remark}

\subsection{Elementary sums}

For \((h,x)\in\mathsf E_{m,d}\), define
\begin{equation}\label{eq:In-def}
  I_t^n(h,x)
  :=
  \sum_{k\ge0}
  h(\tau_k^n(h))
  \Bigl(
    x(t\wedge\tau_{k+1}^n(h))
    -
    x(t\wedge\tau_k^n(h))
  \Bigr).
\end{equation}
By \Cref{lem:crossing-times}, the sum is pathwise finite.  Since \(x\) is
continuous, \(I^n(h,x)\) belongs to \(C([0,T];\R^m)\).

Set \(h(0-)=h(0)\) and define the elementary left-point process
\begin{equation}\label{eq:Hn-def}
  H_s^n(h)
  =
  h(0)\one_{\{0\}}(s)
  +
  \sum_{k\ge0}
  h(\tau_k^n(h))
  \one_{(\tau_k^n(h),\,\tau_{k+1}^n(h)]}(s).
\end{equation}

\begin{lemma}[Uniform approximation of \(h_-\)]\label{lem:integrand-approx}
For every \(h\in\mathsf H_{m,d}\),
\begin{equation}\label{eq:Hn-error}
  \sup_{0\le s\le T}
  \norm{H_s^n(h)-h(s-)}
  \le 2^{-n}.
\end{equation}
If \(P\) is a probability law under which the coordinate \(x\) is a continuous
semimartingale in the usual augmentation of the joint canonical filtration,
then
\begin{equation}\label{eq:In-simple-integral}
  I_t^n(h,x)=\int_0^t H_s^n(h)\,\dd x_s,
  \qquad 0\le t\le T,
\end{equation}
\(P\)-almost surely.
\end{lemma}

\begin{proof}
At \(s=0\), the assertion follows from the convention \(h(0-)=h(0)\).
For \(s>0\), choose the unique active interval
\[
  s\in(\tau_k^n,\tau_{k+1}^n].
\]
For every \(r\in(\tau_k^n,s)\), the definition of
\(\tau_{k+1}^n\) gives
\[
  \norm{h(r)-h(\tau_k^n)}\le\eps_n.
\]
Letting \(r\uparrow s\) and using the existence of \(h(s-)\), we obtain
\[
  \norm{h(s-)-h(\tau_k^n)}\le\eps_n.
\]
This includes \(s=\tau_{k+1}^n\), even when the crossing is caused by a
jump, because the old anchor is used at the right endpoint.  It also includes
\(s=T\), and proves \eqref{eq:Hn-error}.

Since \(h(\tau_k^n)\) is measurable with respect to the stopping-time
sigma-field associated with \((\cF_{t+}^0)\), each process
\[
  h(\tau_k^n)\one_{(\tau_k^n,\tau_{k+1}^n]}
\]
is predictable for the usual augmentation of \((\cF_{t+}^0)\).  To handle
the random number of intervals and the random coefficients, let
\[
  \zeta_N
  =\inf\left\{
    t\in[0,T]:\sup_{0\le u\le t}\norm{h(u)}>N
  \right\}\wedge T.
\]
Truncate \eqref{eq:Hn-def} after \(K\) intervals and stop it at
\(\tau_{K+1}^n\wedge\zeta_N\).  For this integral we may set
\(H_0^n=0\), since changing the integrand at time zero does not affect
integration against the continuous process \(x\).  The resulting integrand is bounded and
elementary predictable, so its stochastic integral is the correspondingly
truncated and stopped sum in \eqref{eq:In-def}.  These identities hold
simultaneously for all \(K,N\) outside one null set.  By
\Cref{lem:crossing-times}, every path has \(\tau_K^n=T\) for all sufficiently
large \(K\); every c\`adl\`ag path is also bounded on \([0,T]\), so
\(\zeta_N=T\) for all sufficiently large \(N\).  Letting first
\(K\to\infty\) and then \(N\to\infty\), and using locality of the
semimartingale integral, gives \eqref{eq:In-simple-integral}.  Continuity of
\(x\) removes any endpoint ambiguity.
\end{proof}

\subsection{A Borel causal convergence system}

For \(t\in[0,T]\), let
\[
  r_t:C([0,T];\R^q)\longrightarrow C([0,t];\R^q)
\]
denote restriction; the same notation is used for c\`adl\`ag paths.  Define
\(\tau_k^{n,t}\) by \eqref{eq:tau-zero}--\eqref{eq:tau-recursion} on
\([0,t]\), with terminal time \(t\), and let \(I^{n,t}\) be the corresponding
sum.  Then
\begin{equation}\label{eq:horizon-consistency}
  r_t I^{n,T}(h,x)
  =
  I^{n,t}(r_th,r_tx).
\end{equation}
Indeed, induction shows that all crossings strictly before \(t\) agree.  If
the next global crossing is greater than or equal to \(t\), both sums use the
same last anchor up to \(t\).  This remains true when an open-passage
infimum equals \(t\) although the threshold is exceeded only after \(t\):
the new interval has zero increment at the endpoint.

Define
\begin{equation}\label{eq:Dt-def}
  \cD_t^{m,d}
  :=
  \left\{
    (h,x)\in\mathsf E_{m,d}:
    (r_tI^n(h,x))_{n\ge1}
    \text{ converges uniformly on }[0,t]
  \right\}.
\end{equation}
When the dimensions are clear, we write simply \(\cD_t\).

\begin{theorem}[Borel and causal convergence domains]
\label{thm:borel-domain}
For every \(t\in[0,T]\), the following statements hold.

\begin{enumerate}[label=(\roman*),leftmargin=2.2em]
  \item The map
  \[
    I^n:\mathsf E_{m,d}\longrightarrow C([0,T];\R^m)
  \]
  is Borel.

  \item The set \(\cD_t\) belongs to \(\cF_t^0\) and is Borel.  Explicitly,
  \begin{equation}\label{eq:Dt-cauchy}
    \cD_t
    =
    \bigcap_{\ell\ge1}
    \bigcup_{N\ge1}
    \bigcap_{n,r\ge N}
    \left\{
      \norm{r_tI^n-r_tI^r}_{\infty}<\ell^{-1}
    \right\}.
  \end{equation}

  \item The partial limit
  \begin{equation}\label{eq:It-partial}
    \cI_t(h,x)
    :=
    \lim_{n\to\infty}r_tI^n(h,x),
    \qquad (h,x)\in\cD_t,
  \end{equation}
  is a Borel map from \(\cD_t\) to \(C([0,t];\R^m)\).

  \item If \(0\le s\le t\le T\), then
  \(\cD_t\subset\cD_s\) and
  \begin{equation}\label{eq:projective-integral}
    r_s\cI_t=\cI_s\qquad\text{on }\cD_t.
  \end{equation}

  \item The system is causal.  If two pairs \((h,x)\) and
  \((\widetilde h,\widetilde x)\) agree on \([0,t]\), then either both belong
  to \(\cD_t\) or neither does; in the former case their \(\cI_t\)-values
  agree.
\end{enumerate}
\end{theorem}

\begin{proof}
By \Cref{lem:crossing-times}, every anchor time is Borel.  Joint Borel
measurability of evaluation shows that
\[
  (h,x)\longmapsto
  h(\tau_k^n(h))
  \bigl(
    x(\,\cdot\wedge\tau_{k+1}^n(h))
    -
    x(\,\cdot\wedge\tau_k^n(h))
  \bigr)
\]
is a Borel \(C([0,T];\R^m)\)-valued map.  The sum stabilizes after finitely
many terms for each path.  Equivalently, it is the pointwise limit of its
finite partial sums, and hence \(I^n\) is Borel.

Since \(C([0,t];\R^m)\) is complete, convergence is equivalent to the Cauchy
condition \eqref{eq:Dt-cauchy}.  Identity
\eqref{eq:horizon-consistency} shows that \(r_tI^{n,T}\) factors through the
Borel restriction map \((h,x)\mapsto(r_th,r_tx)\).  It is therefore
\(\cF_t^0\)-measurable as a \(C([0,t];\R^m)\)-valued map, and
\eqref{eq:Dt-cauchy} gives \(\cD_t\in\cF_t^0\).  The limit of a convergent
sequence of Borel maps into a Polish space is Borel on the convergence set,
which proves (iii).

The horizon-consistency identity proves (iv) and (v) for every elementary
sum.  Passing to the uniform limit proves them for \(\cI_t\).
\end{proof}

\begin{remark}[Factorization through finite horizons]
\label{rem:horizon-factorization}
For each \(t\), the construction used in
\eqref{eq:horizon-consistency} lives directly on
\[
  D([0,t];\R^{m\times d})\times C([0,t];\R^d)
\]
with terminal horizon \(t\).  Denote its domain and limit functional by
\(\widehat{\cD}_t\) and \(\widehat{\cI}_t\).  Identity
\eqref{eq:horizon-consistency} gives
\[
  \cD_t=(r_t\times r_t)^{-1}(\widehat{\cD}_t),
  \qquad
  \cI_t=\widehat{\cI}_t\circ(r_t\times r_t)
  \quad\text{on }\cD_t .
\]
Thus every prefix object below is defined on the finite-horizon
path spaces.  We suppress the hats when no confusion can arise.
\end{remark}

\begin{definition}[Intrinsic It\^o domain and functional]
\label{def:intrinsic-integral}
We set
\[
  \cD^{m,d}:=\cD_T^{m,d},
  \qquad
  \cI:=\cI_T.
\]
Thus
\[
  \cI:\cD^{m,d}\longrightarrow C([0,T];\R^m)
\]
is the endogenous-partition It\^o functional.
\end{definition}

\begin{remark}[Partial causality]
\label{rem:partial-causal}
Extending \(\cI\) by zero on \((\cD^{m,d})^c\) gives a total Borel map.
That extension need not be causal, because membership in \(\cD_T\) uses the
whole time interval.  The intrinsic causal object is the compatible family of
partial maps
\[
  \bigl(\cD_t^{m,d},\cI_t\bigr)_{0\le t\le T}.
\]
\end{remark}

\subsection{Identification under every semimartingale law}

Let \(\mathfrak P_{\mathrm{csm}}^{m,d}\) denote the class of Borel probability
measures \(P\) on \(\mathsf E_{m,d}\) such that the coordinate \(x\) is a
continuous semimartingale relative to the usual augmentation of the joint
canonical filtration generated by \((h,x)\).  This joint-filtration
condition is essential: arbitrary anticipative information carried by \(h\)
can destroy the semimartingale property of \(x\).

\begin{theorem}[Law-free semimartingale identification]
\label{thm:semimartingale-identification}
For every \(P\in\mathfrak P_{\mathrm{csm}}^{m,d}\),
\begin{equation}\label{eq:D-full}
  P(\cD^{m,d})=1
\end{equation}
and
\begin{equation}\label{eq:I-identification}
  \cI(h,x)
  =
  \int_0^\cdot h(s-)\,\dd x_s
  \qquad
  P\text{-almost surely in }C([0,T];\R^m).
\end{equation}
In particular, the full deterministic sequence \(I^n\) converges
\(P\)-almost surely; no \(P\)-dependent subsequence is used.
\end{theorem}

\begin{proof}
Work first after localization.  Write the continuous semimartingale as
\[
  x=M+A,
\]
where \(M\) is a continuous local martingale and \(A\) is a continuous
finite-variation process.  Set
\[
  V_t
  =
  \sum_{j=1}^d[M^j]_t
  +
  \sum_{j=1}^d\Var(A^j;[0,t]),
  \qquad
  \rho_R
  =
  \inf\{t:V_t>R\}\wedge T .
\]
The process \(V\) is continuous, so
\[
  V_{\rho_R}\le R,
  \qquad
  \rho_R\uparrow T
\]
almost surely.

Let \(K^n=H^n-h_-\).  By \Cref{lem:integrand-approx},
\(\norm{K^n}_\infty\le2^{-n}\).  The finite-variation part satisfies
\begin{equation}\label{eq:fv-error}
  \sup_{t\le T}
  \abs{
    \int_0^{t\wedge\rho_R}K_s^n\,\dd A_s
  }
  \le R2^{-n}.
\end{equation}
For the local-martingale part, the Burkholder--Davis--Gundy inequality gives,
for every \(\eta>0\),
\begin{align}
  P\left(
    \sup_{t\le T}
    \abs{\int_0^{t\wedge\rho_R}K_s^n\,\dd M_s}>\eta
  \right)
  &\le
  \eta^{-2}
  \E\left[
    \sup_{t\le T}
    \abs{\int_0^{t\wedge\rho_R}K_s^n\,\dd M_s}^{\,2}
  \right] \notag\\
  &\le
  C_{d,m}\eta^{-2}
  \E\int_0^{\rho_R}\norm{K_s^n}^2\,\dd[M]_s
  \le C_{d,m}\eta^{-2}R\,4^{-n}.
  \label{eq:martingale-error}
\end{align}
The right-hand side is summable in \(n\).  Apply Borel--Cantelli with
\(\eta=j^{-1}\), successively for \(j\in\N\), and intersect the resulting
probability-one events.  Together with \eqref{eq:fv-error}, this yields
\[
  \sup_{t\le T}
  \abs{
    I_{t\wedge\rho_R}^n(h,x)
    -
    \int_0^{t\wedge\rho_R}h(s-)\,\dd x_s
  }
  \longrightarrow0
\]
almost surely along the full sequence for each fixed \(R\).  Intersecting
also over \(R\in\N\) preserves probability one.  For every sample point in
this event, \(V_T<\infty\), so some integer \(R\) satisfies
\(\rho_R=T\).  The stopped convergence is then uniform on \([0,T]\), proving
\eqref{eq:D-full} and \eqref{eq:I-identification}.
\end{proof}

\begin{corollary}[A common polar exceptional set]\label{cor:fixed-polar-integral}
Let
\[
  \cP\subset\mathfrak P_{\mathrm{csm}}^{m,d}
\]
be arbitrary; no domination, compactness, convexity, or mutual absolute
continuity is required.  Then
\[
  N_{\mathrm{int}}
  :=
  \mathsf E_{m,d}\setminus\cD^{m,d}
\]
is a Borel \(\cP\)-polar set.  On its complement, the same Borel path
functional \(\cI\) represents the stochastic integral under every
\(P\in\cP\).
\end{corollary}

\begin{proof}
\Cref{thm:semimartingale-identification} gives
\(P(N_{\mathrm{int}})=0\) for every \(P\in\cP\).  Apply
\Cref{prop:fixed-event}.
\end{proof}

\begin{corollary}[Quasi-sure scheme independence]
\label{cor:scheme-independence}
Suppose \(\widetilde{\cI}\) is a second Borel path functional, defined on a
Borel domain \(\widetilde{\cD}\subset\mathsf E_{m,d}\), such that for every
\(P\in\cP\),
\[
  P(\widetilde{\cD})=1,
  \qquad
  \widetilde{\cI}
  =
  \int_0^\cdot h(s-)\,\dd x_s
  \quad P\text{-almost surely}.
\]
Then
\[
  \bigl\{(h,x)\in\cD\cap\widetilde{\cD}:
    \cI(h,x)\ne\widetilde{\cI}(h,x)\bigr\}
\]
is \(\cP\)-polar.
\end{corollary}

\begin{proof}
Under each \(P\in\cP\), both functionals are versions of the same continuous
stochastic integral.  Their inequality event is therefore \(P\)-null for
every \(P\).
\end{proof}

\begin{remark}[The domain is a pair domain]\label{rem:pair-domain}
The domain in \Cref{cor:fixed-polar-integral} lies in the joint space of
\((h,x)\).  A driver-only domain uniform over an uncountable integrand class
requires additional separability, stability, or capacity estimates, since
projection need not preserve polarity.
\end{remark}

\section{Law-free equations and universal solution maps}
\label{sec:equations}

\subsection{A scheme-fixed quadratic-covariation functional}

The endogenous integral can be used to construct quadratic covariation from
the first-level path alone.  For \(1\le i,j\le d\), let \(e_j\) be the
\(j\)-th coordinate vector in \(\R^d\), and define the row-valued continuous
path
\begin{equation}\label{eq:coordinate-integrand}
  \chi^{ij}(x)_t=x_t^i e_j^\top .
\end{equation}
Set
\begin{equation}\label{eq:Dqv}
  \cD_{\mathrm{qv}}
  =
  \bigcap_{1\le i,j\le d}
  \left\{
    x\in\mathsf X_d:
    \bigl(\chi^{ij}(x),x\bigr)\in\cD^{1,d}
  \right\}.
\end{equation}
For \(x\in\cD_{\mathrm{qv}}\), define
\begin{equation}\label{eq:pathwise-qv}
  \qv_t^{ij}(x)
  =
  x_t^ix_t^j-x_0^ix_0^j
  -
  \cI_t\bigl(\chi^{ij}(x),x\bigr)
  -
  \cI_t\bigl(\chi^{ji}(x),x\bigr).
\end{equation}
Thus \(\qv(x)\) is symmetric by construction.

Let \(C_{\mathrm{fv}}([0,T];\R^{d\times d})\) be the set of continuous
matrix-valued paths whose entries have finite total variation, equipped with
the trace Borel structure inherited from
\(C([0,T];\R^{d\times d})\).  Let \(\Pi_{\mathbb Q}\) be the countable family
of partitions containing \(0,T\) and having rational interior points.  For a
continuous vector path \(a\),
\[
  \Var(a;[0,T])
  =
  \sup_{\pi\in\Pi_{\mathbb Q}}
  \sum_{\ell}
  \norm{a(t_{\ell+1})-a(t_\ell)}.
\]
Thus finite variation is a Borel condition, being expressed through a
countable supremum of continuous evaluation functionals.

\begin{definition}[Scheme-fixed covariation domain]
\label{def:Gqv}
Define
\begin{equation}\label{eq:Gqv}
  \cG_{\mathrm{qv}}
  =
  \left\{
    x\in\cD_{\mathrm{qv}}:
    \qv(x)\in C_{\mathrm{fv}}([0,T];\R^{d\times d})
  \right\}.
\end{equation}
\end{definition}

\begin{theorem}[Pathwise covariation under all semimartingale laws]
\label{thm:pathwise-qv}
The set \(\cG_{\mathrm{qv}}\) is Borel, and
\[
  \qv:\cG_{\mathrm{qv}}
  \longrightarrow C_{\mathrm{fv}}([0,T];\R^{d\times d})
\]
is Borel.  If \(P\) is any law on \(\mathsf X_d\) under which the coordinate
process \(X\) is a continuous semimartingale, then
\begin{equation}\label{eq:qv-identification}
  P(\cG_{\mathrm{qv}})=1,
  \qquad
  \qv^{ij}(X)=[X^i,X^j]
  \quad P\text{-almost surely}
\end{equation}
for every \(i,j\).  Consequently, for every family \(\cP\) of such laws,
\(\mathsf X_d\setminus\cG_{\mathrm{qv}}\) is one Borel
\(\cP\)-polar set.
\end{theorem}

\begin{proof}
The map \(x\mapsto\chi^{ij}(x)\) is continuous.  Hence
\eqref{eq:Dqv}, \Cref{thm:borel-domain}, and the Borel nature of finite
variation show that \(\cG_{\mathrm{qv}}\) and \(\qv\) are Borel.

Under a continuous-semimartingale law, the joint filtration generated by
\((\chi^{ij}(X),X)\) is the canonical filtration of \(X\).  Applying
\Cref{thm:semimartingale-identification} and the integration-by-parts formula
gives
\begin{align*}
  X_t^iX_t^j-X_0^iX_0^j
  &=
  \int_0^tX_s^i\,\dd X_s^j
  +
  \int_0^tX_s^j\,\dd X_s^i
  +
  [X^i,X^j]_t .
\end{align*}
This is exactly \eqref{eq:qv-identification}.  Continuous-semimartingale
covariations have finite variation, so \(P(\cG_{\mathrm{qv}})=1\).
The polar assertion follows from \Cref{prop:fixed-event}.
\end{proof}

\begin{remark}[Dependence on the deterministic scheme]
\label{rem:qv-scheme}
The functional \(\qv\) is fixed by the particular level-crossing sequence
used in \eqref{eq:tau-recursion}; on exceptional deterministic paths another
scheme may give a different value or no value.  Under any family of
continuous-semimartingale laws, two Borel schemes representing classical
covariation agree quasi surely by \Cref{cor:scheme-independence}.  The
scheme-fixed qualifier is therefore essential pathwise and disappears only
after quasi-sure identification.
\end{remark}

\begin{remark}
In dimension \(d>1\), the two stochastic-integral terms in
\eqref{eq:pathwise-qv} are generally different.  The expression
\[
  X_tX_t^\top-X_0X_0^\top
  -2\int_0^tX_s\,\dd X_s^\top
\]
is therefore not the matrix quadratic covariation unless additional
symmetry is present.
\end{remark}

\subsection{Stieltjes integration}

The second deterministic ingredient is a Borel integration map for
finite-variation paths.

\begin{lemma}[Borel Riemann--Stieltjes map]\label{lem:stieltjes-borel}
Let \(p,q\ge1\).  On
\[
  C([0,T];\R^{p\times q})
  \times
  C_{\mathrm{fv}}([0,T];\R^q),
\]
the map
\[
  (f,a)\longmapsto
  \left(t\longmapsto\int_0^t f_s\,\dd a_s\right)
  \in C([0,T];\R^p)
\]
is Borel.
\end{lemma}

\begin{proof}
For the dyadic partition
\(\pi_n=\{t_k=kT2^{-n}:0\le k\le2^n\}\), define
\[
  S_t^n(f,a)
  =
  \sum_{k=0}^{2^n-1}
  f(t_k)
  \bigl[
    a(t\wedge t_{k+1})-a(t\wedge t_k)
  \bigr].
\]
The map \((f,a)\mapsto S^n(f,a)\) is continuous from the product of the two
ambient continuous-path spaces to \(C([0,T];\R^p)\).  For each fixed
continuous \(f\) and finite-variation \(a\), the sums converge uniformly in
\(t\), with
\[
  \norm{
    S^n(f,a)-\int_0^\cdot f_s\,\dd a_s
  }_\infty
  \le
  \omega_f(T2^{-n})\,\Var(a;[0,T]),
\]
where \(\omega_f\) is the modulus of continuity of \(f\).  The limit is
therefore Borel on the stated Borel domain.
\end{proof}

\subsection{A Borel solution relation}

Fix a solution dimension \(q\ge1\), an initial point \(y_0\in\R^q\), and
continuous, locally bounded coefficients
\begin{align}
  b&:[0,T]\times\R^q\longrightarrow\R^q, \notag\\
  \beta^{ij}&:[0,T]\times\R^q\longrightarrow\R^q,
       \qquad 1\le i,j\le d, \notag\\
  \sigma&:[0,T]\times\R^q\longrightarrow\R^{q\times d}.
  \label{eq:coefficients}
\end{align}
For \(y\in\mathsf X_q\), write
\[
  \sigma[y]_t=\sigma(t,y_t).
\]

\begin{lemma}[Prefix covariation system]\label{lem:prefix-qv}
For \(u\in[0,T]\), set
\[
  \mathsf X_d^{(u)}=C([0,u];\R^d),
  \qquad
  \chi_u^{ij}(x)_s=x_s^ie_j^\top,\quad s\le u,
\]
and define
\[
  \cD_{\mathrm{qv}}^{(u)}
  =
  \bigcap_{i,j=1}^d
  \left\{
    x\in\mathsf X_d^{(u)}:
    (\chi_u^{ij}(x),x)\in\widehat{\cD}_u^{1,d}
  \right\}.
\]
For \(x\in\cD_{\mathrm{qv}}^{(u)}\) and \(s\le u\), let
\begin{align}
  \qv_s^{ij,(u)}(x)
  ={}&
  x_s^ix_s^j-x_0^ix_0^j \notag\\
  &-
  \bigl(\widehat{\cI}_u(\chi_u^{ij}(x),x)\bigr)_s
  -
  \bigl(\widehat{\cI}_u(\chi_u^{ji}(x),x)\bigr)_s
  \label{eq:prefix-qv}
\end{align}
and
\[
  \cG_{\mathrm{qv}}^{(u)}
  =
  \left\{
    x\in\cD_{\mathrm{qv}}^{(u)}:
    \qv^{(u)}(x)\in
    C_{\mathrm{fv}}([0,u];\R^{d\times d})
  \right\}.
\]
These domains and maps are Borel.  If \(0\le v\le u\), then
\begin{equation}\label{eq:prefix-qv-consistency}
  r_v\qv^{(u)}(x)=\qv^{(v)}(r_vx)
  \qquad
  \text{for }x\in\cD_{\mathrm{qv}}^{(u)}.
\end{equation}
\end{lemma}

\begin{proof}
Borel measurability follows from
\Cref{rem:horizon-factorization} exactly as in
\Cref{thm:pathwise-qv}.  Projective consistency of the finite-horizon
integrals in \eqref{eq:projective-integral}, applied to both coordinate
integrands in \eqref{eq:prefix-qv}, gives
\eqref{eq:prefix-qv-consistency}.
\end{proof}

\begin{definition}[Law-free It\^o solution relation]
\label{def:solution-relation}
The relation
\[
  \cR_T\subset\mathsf X_d\times\mathsf X_q
\]
consists of all pairs \((x,y)\) such that
\[
  x\in\cG_{\mathrm{qv}},
  \qquad
  (\sigma[y],x)\in\cD^{q,d},
\]
and, for every \(t\in[0,T]\),
\begin{align}
  y_t
  ={}&y_0
  +\int_0^t b(s,y_s)\,\dd s
  +\sum_{i,j=1}^d
      \int_0^t\beta^{ij}(s,y_s)\,\dd\qv_s^{ij}(x)
  \notag\\
  &+\cI_t(\sigma[y],x).
  \label{eq:law-free-equation}
\end{align}
For \(t<T\), the prefix relation \(\cR_t\) is defined on the finite-horizon
spaces by the same conditions on \([0,t]\), using the factorized objects from
\Cref{rem:horizon-factorization}, the domain
\(\cG_{\mathrm{qv}}^{(t)}\), and the covariation
\(\qv^{(t)}\) from \Cref{lem:prefix-qv}.
\end{definition}

Every symbol in \eqref{eq:law-free-equation} is defined from \((x,y)\)
without reference to a probability measure.  Under a continuous
semimartingale law, the equation becomes
\begin{equation}\label{eq:classical-equation}
  dY_t
  =
  b(t,Y_t)\,\dd t
  +\sum_{i,j=1}^d\beta^{ij}(t,Y_t)\,\dd[X^i,X^j]_t
  +\sigma(t,Y_t)\,\dd X_t.
\end{equation}

\begin{proposition}[Borel graph]\label{prop:borel-graph}
For every \(t\in[0,T]\), the relation
\[
  \cR_t\subset
  C([0,t];\R^d)\times C([0,t];\R^q)
\]
is Borel.  The family is compatible with restriction:
\begin{equation}\label{eq:relation-restriction}
  (x,y)\in\cR_T
  \quad\Longrightarrow\quad
  (r_tx,r_ty)\in\cR_t .
\end{equation}
\end{proposition}

\begin{proof}
The coefficient maps
\[
  y\mapsto b(\,\cdot,y_\cdot),\qquad
  y\mapsto\beta^{ij}(\,\cdot,y_\cdot),\qquad
  y\mapsto\sigma[y]
\]
are continuous for the uniform topology.  The Lebesgue integral is
continuous on continuous integrands, the finite-variation integral is Borel
by \Cref{lem:stieltjes-borel}, and the partial It\^o functional is Borel by
\Cref{thm:borel-domain}.  Thus the right-hand side of
\eqref{eq:law-free-equation} is Borel on its Borel domain.  Equality of two
continuous paths is a closed condition, proving that \(\cR_t\) is Borel.
Restriction compatibility follows from
\eqref{eq:projective-integral} and the corresponding locality of Lebesgue and
Riemann--Stieltjes integration.
\end{proof}

\begin{proposition}[Classical solutions enter the same relation]
\label{prop:classical-enters-relation}
Let \(P\) be a continuous-semimartingale law on \(\mathsf X_d\), and suppose
that \(Y^P\) is a continuous strong canonical solution of
\eqref{eq:classical-equation}.  Here this means that \(Y^P\) is an
\(\mathsf X_q\)-valued, \(\cF_T^P\)-measurable random variable and that
\(Y_t^P\) is adapted to the usual augmentation
\((\cF_t^P)_{t\le T}\) of the canonical filtration of \(X\).  Then
\begin{equation}\label{eq:classical-in-relation}
  (X,Y^P)\in\cR_T
  \qquad P\text{-almost surely}.
\end{equation}
\end{proposition}

\begin{proof}
\Cref{thm:pathwise-qv} identifies \(\qv(X)\) with \([X]\) outside a
\(P\)-null set.  Since \(\mathsf X_q\) is a standard Borel space and
\(\cF_T^P\) is the \(P\)-completion of the terminal raw canonical
sigma-field, \(Y^P\) has an \(\cF_T^0\)-Borel
\(\mathsf X_q\)-valued version.  We replace \(Y^P\) by this version, retaining
adaptedness and the equation outside a \(P\)-null set.

Let
\[
  \overline P
  =
  P\circ(\sigma[Y^P],X)^{-1}
\]
be the resulting Borel probability law on
\(\mathsf H_{q,d}\times\mathsf X_d\).  Let \((\mathcal G_t)_{t\le T}\) be
the usual augmentation of the filtration generated by
\((\sigma[Y^P],X)\) on the original canonical space.  Adaptedness of \(Y^P\)
gives \(\mathcal G_t\subset\cF_t^P\), and \(X\) is
\((\mathcal G_t)\)-adapted.  Stricker's theorem \cite{Protter05} implies
that \(X\) remains a continuous semimartingale in \((\mathcal G_t)\).
Equivalently,
\(\overline P\in\mathfrak P_{\mathrm{csm}}^{q,d}\).  The continuous process
\(\sigma[Y^P]\) is predictable in both filtrations, and its stochastic
integrals agree because both are limits of the same elementary predictable
integrals.  Therefore \Cref{thm:semimartingale-identification} identifies
\(\cI\) with the last It\^o integral in \eqref{eq:classical-equation}.
Substitution gives \eqref{eq:law-free-equation} outside one \(P\)-null set.
\end{proof}

\subsection{Selection, aggregation, and causality}

The following selection results are consequences of Jankov--von Neumann and
Lusin--Novikov uniformization.

\begin{theorem}[Analytic existence and universal selection]
\label{thm:analytic-selection}
Let \(E,F\) be Polish spaces, let
\(\cR\subset E\times F\) be analytic, and let \(\cP\) be a family of Borel
probability measures on \(E\).  Suppose that for every \(P\in\cP\) there is
an \(F\)-valued, \(P\)-completion-measurable random variable \(Y^P\) such
that
\[
  (e,Y^P(e))\in\cR
  \qquad P\text{-almost surely}.
\]
Then the intrinsic existence domain
\begin{equation}\label{eq:intrinsic-existence-domain}
  G_{\cR}:=\pi_E(\cR)=\{e:\cR_e\ne\varnothing\}
\end{equation}
is analytic, universally measurable, and satisfies
\[
  P(G_{\cR})=1\qquad\text{for every }P\in\cP,
  \qquad
  c_{\cP}(E\setminus G_{\cR})=0.
\]
There exists a selector
\[
  S:G_{\cR}\longrightarrow F,
  \qquad
  (e,S(e))\in\cR
\]
for every \(e\in G_{\cR}\), measurable from the universal sigma-field on
\(G_{\cR}\) to \(\cB(F)\).
\end{theorem}

\begin{proof}
The projection of an analytic set is analytic, hence universally measurable.
The law-specific solution gives \(P(G_{\cR})=1\).  The
Jankov--von Neumann uniformization theorem provides a universally measurable
selector on the projection; see \cite{Kechris95}.
\end{proof}

\begin{remark}\label{rem:selection-not-causal}
Jankov--von Neumann yields a universally measurable selector.  Borel
measurability, causality, and identification with prescribed measurewise
solutions require additional deterministic structure.
\end{remark}

\begin{theorem}[Borel aggregation under deterministic uniqueness]
\label{thm:borel-aggregation}
In the setting of \Cref{thm:analytic-selection}, assume in addition that
\(\cR\) is Borel and that
\begin{equation}\label{eq:singleton-sections}
  \abs{\cR_e}\le1\qquad\text{for every }e\in E.
\end{equation}
Then \(G_{\cR}\) is Borel, and the unique map
\[
  S:G_{\cR}\longrightarrow F,
  \qquad
  \Graph(S)=\cR,
\]
is Borel.  The fixed Borel set \(E\setminus G_{\cR}\) is
\(\cP\)-polar.  Moreover,
\begin{equation}\label{eq:aggregate-equality}
  S=Y^P\qquad P\text{-almost surely}
\end{equation}
for every law-specific solution \(Y^P\).
\end{theorem}

\begin{proof}
A Borel subset of a product of standard Borel spaces with countable vertical
sections has Borel projection and Borel uniformization by the
Lusin--Novikov theorem; see \cite{Kechris95}.  Under
\eqref{eq:singleton-sections}, this uniformization is the unique map \(S\).
The existence hypothesis gives \(P(G_{\cR})=1\) for every \(P\), so the
complement is polar.  Finally, both \(S(e)\) and \(Y^P(e)\) belong to the same
singleton section for \(P\)-almost every \(e\), which proves
\eqref{eq:aggregate-equality}.
\end{proof}

\begin{corollary}[Quasi-sure aggregation for the law-free equation]
\label{cor:equation-aggregation}
Let \(\cP\) be any family of continuous-semimartingale laws on
\(\mathsf X_d\).  Assume:
\begin{enumerate}[label=(\roman*),leftmargin=2.2em]
  \item under each \(P\in\cP\), \eqref{eq:classical-equation} admits a
  continuous strong canonical solution \(Y^P\);
  \item the law-free relation \(\cR_T\) in
  \Cref{def:solution-relation} has at most one solution for every
  \(x\in\mathsf X_d\).
\end{enumerate}
Then there is a Borel set \(G\subset\mathsf X_d\) with
\[
  P(G)=1\qquad\text{for every }P\in\cP
\]
and a Borel map \(S:G\to\mathsf X_q\) such that
\[
  (x,S(x))\in\cR_T\quad\text{for every }x\in G,
  \qquad
  S(X)=Y^P\quad P\text{-almost surely}.
\]
The pair \((G,S)\) is unique up to a \(\cP\)-polar set.
\end{corollary}

\begin{proof}
\Cref{prop:borel-graph} gives a Borel relation,
\Cref{prop:classical-enters-relation} gives measurewise existence inside that
same relation, and assumption (ii) gives singleton sections.  Apply
\Cref{thm:borel-aggregation}.
\end{proof}

\begin{remark}[Deterministic uniqueness]\label{rem:uniqueness-required}
Assumption~(ii) in \Cref{cor:equation-aggregation} is uniqueness on each
fixed driver path.  It is stronger than pathwise uniqueness under each
individual law.
\end{remark}

\begin{corollary}[Prefix uniqueness implies causality]
\label{cor:prefix-causality}
Assume the hypotheses of \Cref{cor:equation-aggregation} and, in addition,
that every prefix relation \(\cR_t\) has at most one element in each vertical
section.  If \(x,\widetilde x\in G\) agree on \([0,t]\), then
\begin{equation}\label{eq:causal-S}
  r_tS(x)=r_tS(\widetilde x).
\end{equation}
Thus the aggregate is causal on its intrinsic domain.

Extend \(S\) to all of \(\mathsf X_d\) by a fixed continuous path outside
\(G\).  Let
\[
  \mathcal N^P
  =
  \left\{
    N:N\subset N_0\text{ for some }
    N_0\in\cF_T^0,\ P(N_0)=0
  \right\},
\]
\[
  \cF_t^P
  =
  \bigcap_{u>t}\bigl(\cF_u^0\vee\mathcal N^P\bigr),
  \quad t<T,
  \qquad
  \cF_T^P=\cF_T^0\vee\mathcal N^P,
\]
and
\[
  \cF_t^{\cP}
  :=
  \bigcap_{P\in\cP}\cF_t^P,
\]
so that \((\cF_t^P)\) is the usual augmentation of the raw canonical
filtration.
The extended process \(S_t(X)\) is adapted to
\((\cF_t^{\cP})_{t\le T}\) and, having continuous paths, is progressively
measurable for this universal filtration.
\end{corollary}

\begin{proof}
If \(x\) and \(\widetilde x\) have the same prefix, then
\eqref{eq:relation-restriction} places \(r_tS(x)\) and
\(r_tS(\widetilde x)\) in the same vertical section of \(\cR_t\).
Prefix uniqueness proves \eqref{eq:causal-S}.

For every \(P\in\cP\), \Cref{cor:equation-aggregation} gives
\(S(X)=Y^P\), \(P\)-almost surely.  For every Borel
\(B\subset\R^q\),
\[
  \{S_t(X)\in B\}
  \mathbin{\triangle}
  \{Y_t^P\in B\}
  \subset
  \{S(X)\ne Y^P\}\in\mathcal N^P.
\]
Since \(Y^P_t\) is \(\cF_t^P\)-measurable, \(S_t(X)\) is also
\(\cF_t^P\)-measurable.  This holds for every \(P\), so it is
\(\cF_t^{\cP}\)-measurable.  An adapted continuous process is progressively
measurable.

This gives progressiveness relative to the intersection of the completed
filtrations.
\end{proof}

\begin{theorem}[Canonical approximation criterion for solution maps]
\label{thm:canonical-solution-approx}
Let \(E\) be a canonical input path space with Borel restriction maps
\(r_t^E\), let
\[
  F=C([0,T];\R^q)
\]
with the uniform topology, and let \(\Phi_n:E\to F\) be Borel maps.  Suppose:
\begin{enumerate}[label=(\roman*),leftmargin=2.2em]
  \item the convergence set \(G_\Phi\) from
  \eqref{eq:common-convergence-set} has full measure under every
  \(P\in\cP\);
  \item for every \(e\in G_\Phi\), the limit
  \(\Phi(e)=\lim_n\Phi_n(e)\) belongs to a fixed Borel solution relation
  \(\cR_e\);
  \item each \(\Phi_n\) is prefix-compatible in the precise sense that
  \[
    r_t^Ee=r_t^Ee'
    \quad\Longrightarrow\quad
    r_t\Phi_n(e)=r_t\Phi_n(e')
  \]
  for all \(e,e'\in E\) and \(t\le T\).
\end{enumerate}
Then \(G_\Phi^c\) is one Borel \(\cP\)-polar set, \(\Phi\) is a Borel
solution map on \(G_\Phi\), and \(\Phi\) is causal there.  If a measurewise
solution \(Y^P\) is the almost sure limit of \(\Phi_n\) under every
\(P\in\cP\), then \(\Phi=Y^P\), \(P\)-almost surely, for every \(P\).
\end{theorem}

\begin{proof}
\Cref{thm:common-approx} proves the Borel and polar assertions and identifies
the measurewise limits.  Since restriction
\(r_t:F\to C([0,t];\R^q)\) is continuous, prefix compatibility passes to the
uniform limit on \(G_\Phi\), which gives causality.
\end{proof}

\section{Quantitative driver-only capacity cores and a simultaneous scalar flow}
\label{sec:driver-core}

The Bichteler--Karandikar functional is naturally defined on a pair
\((h,x)\).  We specialize to scalar uniformly nondegenerate autonomous
equations.  The construction combines a uniform characteristic bound, a
continuous causal quadratic-variation approximation, and the Lamperti
transform.

We use the following moment-to-tail estimate.

\begin{lemma}[Moment-to-tail conversion]
\label{lem:moment-tail}
Let \(\{(P_i,Y_i):i\in I\}\) be any collection of probability measures and
nonnegative random variables, possibly on different spaces, and let \(A>0\).
\begin{enumerate}[label=(\roman*),leftmargin=2.2em]
  \item If
  \(\sup_{i\in I}\norm{Y_i}_{L^r(P_i)}\le A\sqrt r\) for every \(r\ge2\),
  then
  \[
    \sup_{i\in I}P_i(Y_i>u)
    \le C\exp\!\left(-c\frac{u^2}{A^2}\right),
    \qquad u\ge0.
  \]
  \item If
  \(\sup_{i\in I}\norm{Y_i}_{L^r(P_i)}\le Ar\) for every \(r\ge2\),
  then
  \[
    \sup_{i\in I}P_i(Y_i>u)
    \le C\exp\!\left(-c\frac{u}{A}\right),
    \qquad u\ge0.
  \]
\end{enumerate}
Here and below \(c,C\in(0,\infty)\) are numerical constants unless their
dependence is displayed.
\end{lemma}

\begin{proof}
Markov's inequality gives, uniformly in \(i\),
\(P_i(Y_i>u)\le(A\sqrt r/u)^r\) in (i) and
\(P_i(Y_i>u)\le(Ar/u)^r\) in (ii).  For \(u\) larger than a fixed multiple of
\(A\), choose an integer \(r\ge2\) comparable respectively to
\(u^2/A^2\) and \(u/A\).  This gives the displayed estimates.  Increasing
\(C\) handles the remaining bounded range of \(u/A\).
\end{proof}

The endogenous construction satisfies the following uniform pair-domain
estimate.  Let \(\overline\cP_\Lambda\) be any family
of laws of a scalar pair \((h,X)\) such that \(X\) is a continuous local
martingale in the joint filtration and
\(\dd[X]_t\le\Lambda\,\dd t\).  Write \(I^n(h,X)\) for the elementary sums
from \eqref{eq:In-def}, and set
\[
  c_{\overline\cP_\Lambda}(A)
  =\sup_{\overline P\in\overline\cP_\Lambda}\overline P(A).
\]

\begin{proposition}[Uniform endogenous pair-domain rate]
\label{prop:uniform-BK-rate}
There is a numerical constant \(C\) such that, for every \(r\ge2\) and
\(n\ge1\),
\begin{equation}\label{eq:BK-increment-rate}
  \sup_{\overline P\in\overline\cP_\Lambda}
  \left\|
    \norm{I^{n+1}(h,X)-I^n(h,X)}_\infty
  \right\|_{L^r(\overline P)}
  \le C\sqrt r\,(\Lambda T)^{1/2}2^{-n}.
\end{equation}
Consequently, for every \(0<\gamma<1\), the fixed Borel functional
\[
  Z_\gamma^{\mathrm{BK}}(h,x)
  =
  \sup_{n\ge1}
  2^{\gamma n}
  \norm{I^{n+1}(h,x)-I^n(h,x)}_\infty
\]
satisfies
\begin{equation}\label{eq:BK-defect-tail}
  c_{\overline\cP_\Lambda}
  \bigl(Z_\gamma^{\mathrm{BK}}>R\bigr)
  \le
  C_\gamma
  \exp\!\left(-c_\gamma\frac{R^2}{\Lambda T}\right),
  \qquad R\ge0.
\end{equation}
On \(\{Z_\gamma^{\mathrm{BK}}\le R\}\), the full sequence converges with the
deterministic rate
\[
  \norm{I^n-\cI}_\infty
  \le C_\gamma R2^{-\gamma n}.
\]
This gives a quantitative core on the pair space; the driver-only distinction
is described in \Cref{rem:pair-domain}.
\end{proposition}

\begin{proof}
Let \(H^n\) be the elementary predictable integrand associated with the
level-crossing partition.  By \Cref{lem:integrand-approx},
\[
  \norm{H^{n+1}-H^n}_\infty
  \le 2^{-(n+1)}+2^{-n}
  =3\cdot2^{-(n+1)}.
\]
Since \(I^n=\int H^n\,\dd X\), the sharp-growth form of the
Burkholder--Davis--Gundy inequality \cite{Ren08} gives, for every \(r\ge2\),
\[
  \left\|\norm{I^{n+1}-I^n}_\infty\right\|_{L^r(\overline P)}
  \le
  C\sqrt r\left\|
    \left(\int_0^T
      \abs{H^{n+1}_s-H^n_s}^2\,\dd[X]_s
    \right)^{1/2}
  \right\|_{L^r(\overline P)}
  \le C\sqrt r\,(\Lambda T)^{1/2}2^{-n}.
\]
This application is justified first after localization and then by Fatou's
lemma; the deterministic bracket bound removes the localization.
For \(R\ge(\Lambda T)^{1/2}\),
\Cref{lem:moment-tail} and a union bound give
\[
  c_{\overline\cP_\Lambda}
  \bigl(Z_\gamma^{\mathrm{BK}}>R\bigr)
  \le
  C\sum_{n\ge1}
  \exp\!\left(
    -c\frac{R^2 2^{2(1-\gamma)n}}{\Lambda T}
  \right)
  \le
  C_\gamma\exp\!\left(-c_\gamma\frac{R^2}{\Lambda T}\right),
\]
while the remaining range follows from the trivial capacity bound after
increasing \(C_\gamma\).  This proves \eqref{eq:BK-defect-tail}.
Telescoping the increments on the
sublevel set gives the stated deterministic convergence rate.
\end{proof}

\subsection{Continuous causal realized variation}

Let
\[
  \Omega=C_0([0,T];\R)
\]
with the uniform norm, and write \(X_t(x)=x_t\).  For
\(n\ge1\), put \(t_k^n=kT2^{-n}\).  If
\(t\in[t_j^n,t_{j+1}^n]\), define
\begin{equation}\label{eq:dyadic-qv}
  Q_t^n(x)
  =
  \sum_{k=0}^{j-1}
    \bigl(x_{t_{k+1}^n}-x_{t_k^n}\bigr)^2
  +\bigl(x_t-x_{t_j^n}\bigr)^2 .
\end{equation}
The last, unfinished increment makes \(Q^n(x)\) continuous and preserves
causality, although \(Q^n(x)\) need not be increasing between two grid
points.

\begin{lemma}[Continuous causal approximation]\label{lem:dyadic-continuous}
For every \(n\), the map
\[
  Q^n:\Omega\longrightarrow C([0,T];\R)
\]
is continuous.  It is causal in the sense that
\(r_tx=r_t\widetilde x\) implies
\(r_tQ^n(x)=r_tQ^n(\widetilde x)\).
\end{lemma}

\begin{proof}
For fixed \(n\), formula \eqref{eq:dyadic-qv} contains only finitely many
evaluations and algebraic operations.  If \(x^r\to x\) uniformly, then the
increments converge uniformly and are uniformly bounded; hence every square
and the finite sum converge uniformly in \(t\).  The causality assertion is
immediate from the same formula.
\end{proof}

First define the fixed Cauchy domain
\[
  \cC_{\mathrm{dyad}}
  =\{x\in\Omega:(Q^n(x))_{n\ge1}
       \text{ is Cauchy in the uniform norm}\},
\]
and write \(\widehat Q(x)\) for its uniform limit there.  Define the
driver-only quadratic-variation domain by
\begin{equation}\label{eq:Gdyad}
  \cG_{\mathrm{dyad}}
  =
  \left\{
    x\in\cC_{\mathrm{dyad}}:
    \widehat Q(x)\text{ is increasing and }\widehat Q_0(x)=0
  \right\},
  \qquad Q=\widehat Q\big|_{\cG_{\mathrm{dyad}}}.
\end{equation}
The Cauchy criterion and the closedness of the cone of increasing continuous
paths show that both domains are Borel.  The maps in
\Cref{lem:dyadic-continuous} also show that \(\widehat Q\) on
\(\cC_{\mathrm{dyad}}\), and hence \(Q\) on \(\cG_{\mathrm{dyad}}\), are
Borel and causal.

Fix \(\Lambda>0\).  Let \(\mathfrak M_\Lambda\) be the maximal family of all
Borel probability measures \(P\) on \(\Omega\) such that, in the usual
augmentation of its natural filtration, the coordinate process is a
continuous local martingale and there is a predictable process \(a^P\) with
\begin{equation}\label{eq:volatility-bound}
  \begin{aligned}
    [X]^P_t&=\int_0^t a_s^P\,\dd s,
      &&t\in[0,T],\quad P\text{-almost surely},\\
    0&\le a_s^P\le\Lambda,
      &&P\otimes\dd s\text{-almost everywhere}.
  \end{aligned}
\end{equation}
Set
\begin{equation}\label{eq:maximal-capacity}
  c_\Lambda(A)=\sup_{P\in\mathfrak M_\Lambda}P(A),
  \qquad A\in\mathcal B(\Omega).
\end{equation}
The uncountable subfamily of constant-volatility laws is already
nondominated.  Every result below remains true, with the same constants, when
\(\mathfrak M_\Lambda\) is replaced by an arbitrary subfamily.

\begin{theorem}[Uniform realized-variation estimate]
\label{thm:uniform-realized-qv}
There is a numerical constant \(C\) such that, for every \(r\ge2\) and
\(n\ge1\),
\begin{equation}\label{eq:qv-Lr-rate}
  \sup_{P\in\mathfrak M_\Lambda}
  \left\|
    \norm{Q^n(X)-[X]^P}_\infty
  \right\|_{L^r(P)}
  \le
  Cr\Lambda T2^{-n/2}.
\end{equation}
In particular, there are numerical constants \(c,C>0\) such that
\begin{equation}\label{eq:qv-exponential-tail}
  \sup_{P\in\mathfrak M_\Lambda}
  P\!\left(
    \norm{Q^n(X)-[X]^P}_\infty>u
  \right)
  \le
  C\exp\!\left(
    -c\frac{2^{n/2}u}{\Lambda T}
  \right),
  \qquad u\ge0.
\end{equation}
Consequently,
\begin{equation}\label{eq:Gdyad-full-capacity}
  P(\cG_{\mathrm{dyad}})=1
  \qquad\text{for every }P\in\mathfrak M_\Lambda.
\end{equation}
Moreover, for each \(P\in\mathfrak M_\Lambda\), the intrinsic limit \(Q(X)\) is a
version of the \(P\)-quadratic variation:
\begin{equation}\label{eq:intrinsic-qv-identification}
  Q(X)=[X]^P
  \qquad P\text{-almost surely},
\end{equation}
where \([X]^P\) denotes any \(P\)-quadratic-variation version.
\end{theorem}

\begin{proof}
Let \(\eta_n(s)=t_k^n\) on
\([t_k^n,t_{k+1}^n)\).  It\^o's formula on each dyadic interval, including
the unfinished interval containing \(t\), gives the exact identity
\begin{equation}\label{eq:qv-error-martingale}
  Q_t^n(X)-[X]^P_t
  =2\int_0^t\bigl(X_s-X_{\eta_n(s)}\bigr)\,\dd X_s .
\end{equation}
The bracket bound makes \(X\) a true martingale with moments of every order.
The standard sharp-growth BDG estimate \cite{Ren08}
\[
  \left\|\sup_{t\le T}\abs{M_t}\right\|_{L^r}
  \le C\sqrt r\,
     \left\|[M]_T^{1/2}\right\|_{L^r},
  \qquad r\ge2,
\]
applied first to \eqref{eq:qv-error-martingale} gives
\begin{align*}
  \left\|\norm{Q^n(X)-[X]^P}_\infty\right\|_{L^r(P)}
  &\le C\sqrt r\,\Lambda^{1/2}
    \left\|
      \left(\int_0^T
        \abs{X_s-X_{\eta_n(s)}}^2\,\dd s
      \right)^{1/2}
    \right\|_{L^r(P)} \\
  &\le C\sqrt r\,\Lambda^{1/2}
    \left(
      \int_0^T
      \norm{X_s-X_{\eta_n(s)}}_{L^r(P)}^2\,\dd s
    \right)^{1/2}.
\end{align*}
A second application, now to each martingale increment, yields
\[
  \norm{X_s-X_{\eta_n(s)}}_{L^r(P)}
  \le C\sqrt r\,
       \bigl(\Lambda(s-\eta_n(s))\bigr)^{1/2}
  \le C\sqrt r\,(\Lambda T2^{-n})^{1/2}.
\]
All estimates may first be applied to stopped martingales; Fatou's lemma and
the deterministic bracket bound then remove the stopping.  Substitution
proves \eqref{eq:qv-Lr-rate}, uniformly in \(P\), and
\Cref{lem:moment-tail}(ii) gives \eqref{eq:qv-exponential-tail}.
Fix \(P\) and one continuous version
\([X]^P\).  Markov's inequality makes
\[
  \sum_{n\ge1}
  P\bigl(\norm{Q^n-[X]^P}_\infty>\eps\bigr)
  <\infty
\]
for every \(\eps>0\).  Borel--Cantelli, followed by a countable intersection
over rational \(\eps\), proves \(Q^n\to[X]^P\) uniformly
\(P\)-almost surely.  Since \([X]^P\) is continuous and increasing, this
proves \eqref{eq:Gdyad-full-capacity}; hence the complement of the fixed
Borel set \(\cG_{\mathrm{dyad}}\) is polar.  The pathwise limit on this
domain is \(Q\), so the same convergence proves
\eqref{eq:intrinsic-qv-identification}.
\end{proof}

\subsection{Two-scale compact driver cores}

We construct compact cores for \(\mathfrak M_\Lambda\).

\begin{theorem}[Two-scale compact driver cores]
\label{thm:driver-core}
Fix
\begin{equation}\label{eq:eta-beta-choice}
  0<\eta<\frac12,
  \qquad
  0<\gamma<\min\left\{\frac12,2\eta\right\}.
\end{equation}
For every \(A,B\ge1\), there is an explicit compact set
\(\cK_{A,B}\subset\cG_{\mathrm{dyad}}\), depending only on
\((\eta,\gamma,\Lambda,T,A,B)\).  The family is increasing in both
parameters:
\[
  A\le A',\quad B\le B'
  \quad\Longrightarrow\quad
  \cK_{A,B}\subseteq\cK_{A',B'}.
\]
There are positive constants depending only on the indicated exponents such
that
\begin{equation}\label{eq:two-scale-capacity}
  c_\Lambda(\cK_{A,B}^c)
  \le C_\eta e^{-c_\eta A^2}+C_\gamma e^{-c_\gamma B}.
\end{equation}
The map
\[
  Q:\cK_{A,B}\longrightarrow C([0,T])
\]
is uniformly continuous, and
\begin{equation}\label{eq:two-scale-qv-rate}
  \sup_{x\in\cK_{A,B}}\norm{Q^n(x)-Q(x)}_\infty
  \le C_\gamma B\Lambda T2^{-\gamma n},
  \qquad n\ge1,
\end{equation}
and, with \(\alpha=\gamma/(1-\eta+\gamma)>0\),
\begin{equation}\label{eq:two-scale-qv-modulus}
  \norm{Q(x)-Q(\widetilde x)}_\infty
  \le C_{\eta,\gamma}A^\alpha B^{1-\alpha}
       (\Lambda T)^{1-\alpha/2}
       \norm{x-\widetilde x}_\infty^\alpha,
  \qquad x,\widetilde x\in\cK_{A,B}.
\end{equation}
Stopping preserves each \(\cK_{A,B}\).  For deterministic concatenation,
there is \(\kappa_\eta<\infty\) such that, for every \(s\in[0,T]\),
\begin{equation}\label{eq:two-scale-concatenation}
  x,\widetilde x\in\cK_{A,B}
  \quad\Longrightarrow\quad
  x\otimes_s\widetilde x
  \in\cK_{2^{1-\eta}A,\,3B+\kappa_\eta A^2},
\end{equation}
where
\[
  (x\otimes_s\widetilde x)_t
  =x_t\one_{\{t\le s\}}
   +\bigl(x_s+\widetilde x_t-\widetilde x_s\bigr)
      \one_{\{t>s\}}.
\]

\end{theorem}

\begin{proof}
Put \(\lambda=\Lambda T\) and write
\[
  \norm{x}_\eta
  =\sup_{0\le s<t\le T}
    \frac{\abs{x_t-x_s}}{\abs{t-s}^\eta}
\]
and introduce the master approximation defect
\begin{equation}\label{eq:master-defect}
  Z_\gamma(x)
  =\sup_{n\ge1}
    2^{\gamma n}\norm{Q^{n+1}(x)-Q^n(x)}_\infty,
\end{equation}
with value \(+\infty\) allowed.  By
\Cref{lem:dyadic-continuous}, the two sublevel sets in the following
definition are closed.  Put
\[
  D_{A,B}
  =\left\{
    \norm{x}_\eta\le A\sqrt\Lambda\,T^{1/2-\eta},\quad
    Z_\gamma(x)\le B\lambda
  \right\}.
\]
On this closed set the sequence \(Q^n\) converges uniformly and
\begin{equation}\label{eq:tail-master-defect}
  \norm{Q^n(x)-\widehat Q(x)}_\infty
  \le
  B\lambda\sum_{m=n}^\infty2^{-\gamma m}
  \le C_\gamma B\lambda2^{-\gamma n}.
\end{equation}
Consequently \(D_{A,B}\subset\cC_{\mathrm{dyad}}\) and
\(\widehat Q:D_{A,B}\to C([0,T])\) is continuous.

Let \(\mathcal A_\Lambda\) be the closed set of increasing paths
\(q\in C([0,T])\) with \(q_0=0\) and
\(\abs{q_t-q_s}\le\Lambda\abs{t-s}\).  Define
\begin{equation}\label{eq:explicit-core}
  \cK_{A,B}
  =D_{A,B}\cap \widehat Q^{-1}(\mathcal A_\Lambda).
\end{equation}
Thus \(\cK_{A,B}\) is closed and contained in \(\cG_{\mathrm{dyad}}\), with
\(Q=\widehat Q\) there.  The H\"older bound and \(x_0=0\) make it compact by
Arzel\`a--Ascoli.  The family is increasing in \((A,B)\).

For the capacity estimate, define the dyadic oscillation functional
\[
  \mathscr H_\eta(x)
  =
  \sup_{m\ge0}\max_{0\le k<2^m}
  \frac{
    \abs{x_{t_{k+1}^m}-x_{t_k^m}}
  }{(T2^{-m})^\eta}.
\]
A deterministic dyadic chaining argument gives
\begin{equation}\label{eq:dyadic-chaining-holder}
  \norm{x}_\eta\le C_\eta\mathscr H_\eta(x).
\end{equation}
Indeed, for fixed \(s<t\), choose the first dyadic scale comparable with
\(t-s\), approximate both endpoints along the nested finer grids, and sum at
most two increments at each level; the resulting geometric series is bounded
by \(C_\eta\mathscr H_\eta(x)\abs{t-s}^\eta\).

The exponential supermartingale inequality and
\([X]^P_t-[X]^P_s\le\Lambda(t-s)\) imply, uniformly in
\(P\in\mathfrak M_\Lambda\),
\begin{equation}\label{eq:martingale-increment-gaussian}
  P\bigl(\abs{X_t-X_s}>u\bigr)
  \le2\exp\!\left(-\frac{u^2}{2\Lambda(t-s)}\right),
  \qquad 0\le s<t\le T.
\end{equation}
Combining \eqref{eq:dyadic-chaining-holder} with a union bound over the
\(2^m\) adjacent increments at level \(m\) gives
\begin{align}
  c_\Lambda\!\left(
    \norm{X}_\eta>A\sqrt\Lambda\,T^{1/2-\eta}
  \right)
  &\le
  2\sum_{m\ge0}2^m
  \exp\!\left(
    -c_\eta A^2 2^{(1-2\eta)m}
  \right) \notag\\
  &\le C_\eta\exp(-c_\eta A^2),
  \qquad A\ge1.
  \label{eq:Holder-capacity}
\end{align}
The last step uses only \(1-2\eta>0\); for \(A\ge1\), the elementary
series is bounded by \(C_\eta e^{-c_\eta A^2}\).

Moreover, \eqref{eq:qv-Lr-rate} and the triangle inequality imply the
linear-in-\(r\) moment bound
\[
  \sup_{P\in\mathfrak M_\Lambda}
  \left\|
    \norm{Q^{n+1}-Q^n}_\infty
  \right\|_{L^r(P)}
  \le Cr\Lambda T2^{-n/2},
  \qquad r\ge2.
\]
Hence \Cref{lem:moment-tail}(ii) and \(\gamma<1/2\) give
\begin{align}
  c_\Lambda(Z_\gamma>B\lambda)
  &\le
  C\sum_{n\ge1}
  \exp\!\left(
    -cB 2^{(1/2-\gamma)n}
  \right) \notag\\
  &\le C_\gamma\exp(-c_\gamma B),
  \qquad B\ge1.
  \label{eq:defect-capacity}
\end{align}
By \eqref{eq:intrinsic-qv-identification}, under each
\(P\in\mathfrak M_\Lambda\) the intrinsic
limit belongs to \(\mathcal A_\Lambda\) almost surely.  Hence the fixed Borel
set
\[
  F_{A,B}
  =\{x\in D_{A,B}:\widehat Q(x)\notin\mathcal A_\Lambda\}
\]
is polar.  Since
\[
  \cK_{A,B}^c
  \subset
  \left\{\norm{x}_\eta>
    A\sqrt\Lambda\,T^{1/2-\eta}\right\}
  \cup\{Z_\gamma>B\lambda\}\cup F_{A,B},
\]
combining
\eqref{eq:Holder-capacity} and \eqref{eq:defect-capacity} proves
\eqref{eq:two-scale-capacity}, while
\eqref{eq:tail-master-defect} proves
\eqref{eq:two-scale-qv-rate}.

If \(x,\widetilde x\in\cK_{A,B}\),
\(e=\norm{x-\widetilde x}_\infty\), and
\(\Delta_n=T2^{-n}\), then the finite-sum formula and the H\"older bounds give
\begin{equation}\label{eq:Qn-modulus}
  \norm{Q^n(x)-Q^n(\widetilde x)}_\infty
  \le C_\eta A\sqrt\lambda\,e\,2^{(1-\eta)n}.
\end{equation}
Indeed, each difference of squares is bounded by the difference of the two
increments, at most \(2e\), times the sum of their sizes, at most
  \(2A\sqrt\Lambda\,T^{1/2-\eta}\Delta_n^\eta\), and there are at most
  \(2^n+1\) terms.  If \(0<e\le(B/A)\sqrt\lambda\), combining
\eqref{eq:Qn-modulus} with the two tails
\eqref{eq:two-scale-qv-rate} and choosing
\(2^n\) comparable to
\(((B/A)\sqrt\lambda/e)^{1/(1-\eta+\gamma)}\) proves
\eqref{eq:two-scale-qv-modulus}.  The case \(e=0\) follows by continuity.
If \(e>(B/A)\sqrt\lambda\), then both \(Q(x)\) and \(Q(\widetilde x)\)
lie in \(\mathcal A_\Lambda\), so their uniform distance is at most
\(\lambda\).  At the threshold \(e=(B/A)\sqrt\lambda\), the right-hand side
of \eqref{eq:two-scale-qv-modulus} equals a constant multiple of
\(B\lambda\), which dominates \(\lambda\) because \(B\ge1\).  This proves
the estimate for all \(e\).

For a stopped path \(x^s_t=x_{t\wedge s}\), formula
\eqref{eq:dyadic-qv} gives exactly
\[
  Q^n(x^s)=Q^n(x)^s.
\]
Hence neither the H\"older seminorm nor \(Z_\gamma\) increases under stopping,
and \(\mathcal A_\Lambda\) is also preserved.  This proves stopping
stability, including pathwise stopping at an arbitrary value of \(s\).

For concatenation, let \(z=x\otimes_s\widetilde x\).  Its H\"older seminorm
is at most \(2^{1-\eta}A\sqrt\Lambda\,T^{1/2-\eta}\): only the case
\(r<s<t\) needs attention, and there
\[
  \abs{z_t-z_r}
  \le A\sqrt\Lambda\,T^{1/2-\eta}
       \bigl((s-r)^\eta+(t-s)^\eta\bigr)
  \le2^{1-\eta}A\sqrt\Lambda\,T^{1/2-\eta}(t-r)^\eta.
\]
Put \(\Delta_n=T2^{-n}\).  If \(s<T\), let
\(a_n\le s<b_n\) be the unique half-open dyadic cell of length
\(\Delta_n\); when \(s\) is a grid point take \(a_n=s\).  Set
\[
  A_n=x_s-x_{a_n},\qquad
  C_n=\widetilde x_s-\widetilde x_{a_n},\qquad
  B_{n,t}
  =\one_{\{t>s\}}
   \bigl(\widetilde x_{t\wedge b_n}-\widetilde x_s\bigr),
\]
and let \(E_t^{n,s}=2(A_n-C_n)B_{n,t}\).  For \(s=T\), set
\(E^{n,T}=0\).  Directly separating the complete cells on the two sides of
\(s\) gives the exact identity
\[
  Q^n(z)
  =Q^n(x^s)+Q^n(\widetilde x)-Q^n(\widetilde x^s)+E^{n,s},
  \qquad
  \norm{E^{n,s}}_\infty
  \le4A^2\lambda\,2^{-2\eta n}.
\]
Indeed, the formula on the sole crossing cell is just the polarization
identity \((u+v)^2-u^2-v^2=2uv\); it also shows that the error vanishes when
\(s\) is a grid point.

Write \(D_n(u)=Q^{n+1}(u)-Q^n(u)\).  Stopping does not increase
\(\norm{D_n(u)}_\infty\), and the preceding identity at levels \(n\) and
\(n+1\) gives
\begin{align*}
  \norm{D_n(z)}_\infty
  &\le
  \norm{D_n(x^s)}_\infty
  +\norm{D_n(\widetilde x)}_\infty
  +\norm{D_n(\widetilde x^s)}_\infty \\
  &\quad
  +\norm{E^{n+1,s}}_\infty+\norm{E^{n,s}}_\infty \\
  &\le
  3B\lambda2^{-\gamma n}
  +4(1+2^{-2\eta})A^2\lambda2^{-2\eta n}.
\end{align*}
It follows from \(\gamma<2\eta\) that
\[
  Z_\gamma(z)
  \le \bigl(3B+\kappa_\eta A^2\bigr)\lambda,
  \qquad
  \kappa_\eta=4(1+2^{-2\eta}).
\]
Passing to the limit gives
\[
  Q_t(z)
  =Q_{t\wedge s}(x)
   +\one_{\{t>s\}}
      \bigl(Q_t(\widetilde x)-Q_s(\widetilde x)\bigr).
\]
For \(r<s<t\), the two increments on the right are bounded by
\(\Lambda(s-r)\) and \(\Lambda(t-s)\); hence this limit remains in
\(\mathcal A_\Lambda\).  This proves
\eqref{eq:two-scale-concatenation} and completes the proof.
\end{proof}

\begin{corollary}[Parabolic diagonal]
\label{cor:parabolic-core}
Define
\begin{equation}\label{eq:parabolic-core}
  \cK_R:=\cK_{R,R^2},
  \qquad R\ge1.
\end{equation}
Then \((\cK_R)_{R\ge1}\) is increasing, and there are
\(c,C,C_*>0\), depending only on \((\eta,\gamma)\), such that
\begin{equation}\label{eq:core-capacity-rate}
  c_\Lambda(\cK_R^c)\le Ce^{-cR^2},
  \qquad R\ge1,
\end{equation}
\begin{equation}\label{eq:qv-core-rate}
  \sup_{x\in\cK_R}\norm{Q^n(x)-Q(x)}_\infty
  \le C_\gamma R^2\Lambda T2^{-\gamma n},
  \qquad n\ge1,
\end{equation}
\begin{equation}\label{eq:qv-uniform-on-core}
  \lim_{n\to\infty}
  \sup_{x\in\cK_R}\norm{Q^n(x)-Q(x)}_\infty=0,
\end{equation}
and
\begin{equation}\label{eq:qv-core-modulus}
  \norm{Q(x)-Q(\widetilde x)}_\infty
  \le C_{\eta,\gamma}R^{2-\alpha}
       (\Lambda T)^{1-\alpha/2}
       \norm{x-\widetilde x}_\infty^\alpha,
  \qquad x,\widetilde x\in\cK_R.
\end{equation}
The diagonal cores are stable under stopping and, for every \(s\in[0,T]\),
\begin{equation}\label{eq:core-concatenation}
  x,\widetilde x\in\cK_R
  \quad\Longrightarrow\quad
  x\otimes_s\widetilde x\in\cK_{C_*R}.
\end{equation}
\end{corollary}

\begin{proof}
Apply \Cref{thm:driver-core} with \((A,B)=(R,R^2)\).  For concatenation,
take
\[
  C_*=\max\left\{2^{1-\eta},(3+\kappa_\eta)^{1/2}\right\}.
\]
\end{proof}

\begin{remark}[Mixed first- and second-order scaling]
The parameters \(A\) and \(B\) retain the natural tail scales of the
first-order H\"older norm and the second-order approximation defect.  The
parabolic choice \(B=A^2\) matches the sub-Gaussian tail in \(A\) with the
subexponential tail in \(B\), yielding \eqref{eq:core-capacity-rate}.
\end{remark}

\begin{remark}[Quasi-continuity]
In the terminology of quasi-continuity, \Cref{thm:driver-core} constructs
compact continuity cores directly from the characteristic bound
\eqref{eq:volatility-bound}, together with rates and stability under stopping
and concatenation; compare \cite{DenisHuPeng11}.
\end{remark}

\subsection{A deterministic Lamperti--F\"ollmer flow}

Let \(\Theta\) be a compact set of pairs
\(\theta=(b_\theta,\sigma_\theta)\) in
\(C_b^1(\R)\times C_b^2(\R)\), equipped with the product norm topology.
Here
\(\norm{f}_{C_b^k}=\sum_{j=0}^k\norm{f^{(j)}}_\infty\).
Assume that, for constants \(L<\infty\) and
\(0<\underline\sigma\le\overline\sigma<\infty\),
\begin{equation}\label{eq:coefficient-class}
  \sup_{\theta\in\Theta}
  \left(
    \norm{b_\theta}_{C_b^1}
    +\norm{\sigma_\theta}_{C_b^2}
  \right)
  \le L,
  \qquad
  \underline\sigma
  \le\sigma_\theta(y)\le\overline\sigma
\end{equation}
for every \((\theta,y)\in\Theta\times\R\).

Define the Lamperti diffeomorphism and its inverse by
\begin{equation}\label{eq:Lamperti}
  F_\theta(y)=\int_0^y\frac{\dd r}{\sigma_\theta(r)},
  \qquad
  G_\theta=F_\theta^{-1},
\end{equation}
and set
\begin{equation}\label{eq:transformed-coefficients}
  \widehat b_\theta(z)
  =\frac{b_\theta(G_\theta(z))}
         {\sigma_\theta(G_\theta(z))},
  \qquad
  \widehat c_\theta(z)
  =-\frac12\sigma_\theta'(G_\theta(z)).
\end{equation}
The functions in \eqref{eq:transformed-coefficients} are bounded and
globally Lipschitz, uniformly in \(\theta\).  In addition,
\[
  G_\theta'(z)=\sigma_\theta(G_\theta(z)),\qquad
  G_\theta''(z)
  =\sigma_\theta'(G_\theta(z))\sigma_\theta(G_\theta(z)),
\]
and \eqref{eq:coefficient-class} implies
\begin{equation}\label{eq:G-third-bound}
  \sup_{\theta\in\Theta}\norm{G_\theta'''}_\infty<\infty .
\end{equation}

Let \(C_\uparrow([0,T])\) be the cone of continuous increasing paths that
start at zero.  For \(x\in\Omega\), \(q\in C_\uparrow([0,T])\),
\(0\le s\le t\le T\), and \(y\in\R\), let
\(Z^{\theta,s,y}(x,q)\) be the unique solution on \([s,T]\) of
\begin{align}
  Z_t^{\theta,s,y}(x,q)
  ={}&F_\theta(y)+x_t-x_s
  +\int_s^t\widehat b_\theta
       \bigl(Z_r^{\theta,s,y}(x,q)\bigr)\,\dd r
  \notag\\
  &+\int_s^t\widehat c_\theta
       \bigl(Z_r^{\theta,s,y}(x,q)\bigr)\,\dd q_r,
  \label{eq:Lamperti-equation}
\end{align}
and define
\begin{equation}\label{eq:Phi-flow}
  \Phi_{s,t}^\theta(x,q;y)
  =G_\theta\bigl(Z_t^{\theta,s,y}(x,q)\bigr).
\end{equation}

\begin{lemma}[Deterministic parameterized flow]
\label{lem:deterministic-flow}
Equation \eqref{eq:Lamperti-equation} has a unique global continuous
solution.  On every set \(\{q_T\le A\}\), the map
\begin{equation}\label{eq:joint-flow-continuity}
  (x,q,\theta,s,t,y)
  \longmapsto
  \Phi_{s,t}^\theta(x,q;y)
\end{equation}
is jointly continuous, locally uniformly in \(y\), when \(x\) and \(q\)
carry the uniform topology.  It is causal and satisfies, for
\(s\le u\le t\),
\begin{equation}\label{eq:cocycle}
  \Phi_{s,t}^\theta(x,q;y)
  =
  \Phi_{u,t}^\theta
  \bigl(x,q;\Phi_{s,u}^\theta(x,q;y)\bigr).
\end{equation}
For fixed \((x,q,\theta,s,t)\), the map
\(y\mapsto\Phi_{s,t}^\theta(x,q;y)\) is a \(C^1\)-diffeomorphism of
\(\R\), and it preserves orientation.  More precisely,
\begin{equation}\label{eq:flow-derivative}
  \partial_y\Phi_{s,t}^\theta(x,q;y)
  =
  \frac{\sigma_\theta(\Phi_{s,t}^\theta(x,q;y))}
       {\sigma_\theta(y)}
  \exp\!\left\{
    \int_s^t\widehat b_\theta'(Z_r)\,\dd r
    +\int_s^t\widehat c_\theta'(Z_r)\,\dd q_r
  \right\}>0.
\end{equation}
Consequently, on \(\{q_T\le A\}\) there are constants
\(0<c_A\le C_A<\infty\), independent of
\((x,\theta,s,t,y)\), such that
\begin{equation}\label{eq:flow-bi-Lipschitz}
  c_A\le\partial_y\Phi_{s,t}^\theta(x,q;y)\le C_A.
\end{equation}
\end{lemma}

\begin{proof}
Existence and uniqueness follow from Picard iteration and Gronwall's
inequality for the finite measure \(\dd r+\dd q_r\).  The uniform bounds in
\eqref{eq:coefficient-class} prevent explosion.

For continuity, suppose that
\((x^n,q^n,\theta_n,s_n,y_n)\) converges to
\((x,q,\theta,s,y)\), with \(q_T^n\le A\).  Uniform convergence of increasing
continuous paths implies weak convergence of the Stieltjes measures
\(\dd q^n\) to \(\dd q\).  We use the following uniform form of
Helly--Bray: if \(f_n\to f\) and \(q^n\to q\) uniformly, the \(q^n\) are
increasing with uniformly bounded terminal values, and \(s_n\to s\), then,
with the usual signed-integral convention when the upper endpoint is below
the lower one,
\begin{equation}\label{eq:uniform-Helly-Bray}
  \sup_{t\in[0,T]}
  \left|
    \int_{s_n}^t f_n(r)\,\dd q_r^n
    -
    \int_s^t f(r)\,\dd q_r
  \right|
  \longrightarrow0.
\end{equation}
Indeed, first replace \(f_n\) by \(f\), at a cost bounded by
\(\norm{f_n-f}_\infty q_T^n\).  Approximate \(f\) uniformly by a \(C^1\)
function \(g\).  Integration by parts expresses
\(\int_0^t g\,\dd q^n\) through \(g(t)q_t^n\) and
\(\int_0^t q^n\,\dd g\), which converge uniformly in \(t\); the uniform
approximation error is bounded by the common total-variation bound.  Finally,
evaluation at \(s_n\to s\) handles the moving lower endpoint.

We also record the coefficient dependence.  Convergence in the product
\(C_b^1\times C_b^2\) norm, together with the common positive lower bound on
\(\sigma\), gives \(F_{\theta_n}\to F_\theta\) locally uniformly.  Since
\(\abs{G_\vartheta(z)}\le\overline\sigma\abs z\) for every
\(\vartheta\in\Theta\), stability of inverses yields
\(G_{\theta_n}\to G_\theta\) locally uniformly.  Substitution in
\eqref{eq:transformed-coefficients} then gives
\[
  \widehat b_{\theta_n}\longrightarrow\widehat b_\theta,
  \qquad
  \widehat c_{\theta_n}\longrightarrow\widehat c_\theta
\]
locally uniformly, with one common Lipschitz constant.

Extend each solution constantly to the left of its starting point by setting
\(Z_t^n=F_{\theta_n}(y_n)\) for \(t\le s_n\).  The integral equations, the
uniform coefficient bounds, uniform convergence of \(x^n,q^n\), and
continuity of the limiting paths make these extensions uniformly bounded and
equicontinuous.  Every uniformly convergent subsequence has a limit solving
\eqref{eq:Lamperti-equation}: the Lebesgue terms pass by uniform convergence,
while the Stieltjes terms pass by \eqref{eq:uniform-Helly-Bray}.  Uniqueness
identifies the limit, so the entire sequence converges uniformly.  The same
sequential argument, with \(y_n\) in any fixed compact interval, proves local
uniformity in \(y\) and hence \eqref{eq:joint-flow-continuity}.

Causality is immediate from the integral equation.  Solving first from
\(s\) to \(u\) and then from \(u\) to \(t\) gives another solution of the
same equation, so uniqueness proves \eqref{eq:cocycle}.

For solutions starting from Lamperti coordinates \(u\) and
\(\widetilde u\), put \(D=Z^u-Z^{\widetilde u}\).  The mean-value theorem
gives bounded continuous coefficients
\[
  a_r=\int_0^1
    \widehat b_\theta'
      \bigl(Z_r^{\widetilde u}+\lambda D_r\bigr)\,\dd\lambda,
  \qquad
  c_r=\int_0^1
    \widehat c_\theta'
      \bigl(Z_r^{\widetilde u}+\lambda D_r\bigr)\,\dd\lambda,
\]
such that
\[
  D_t=D_s+\int_s^t a_rD_r\,\dd r+\int_s^t c_rD_r\,\dd q_r.
\]
The scalar linear equation has the explicit solution
\[
  D_t=D_s
  \exp\!\left\{
    \int_s^t a_r\,\dd r+\int_s^t c_r\,\dd q_r
  \right\}.
\]
If \(K\) bounds \(\abs a+\abs c\), this gives
\begin{equation}\label{eq:two-sided-flow-bound}
  e^{-K(T+A)}\abs{u-\widetilde u}
  \le
  \abs{Z_t^u-Z_t^{\widetilde u}}
  \le
  e^{K(T+A)}\abs{u-\widetilde u}.
\end{equation}
Since \(F_\theta(y)\to\pm\infty\) as \(y\to\pm\infty\), the lower bound in
\eqref{eq:two-sided-flow-bound} shows that the images also tend to
\(\pm\infty\).  Thus the map is an orientation-preserving homeomorphism.
The coefficients \(\widehat b_\theta,\widehat c_\theta\) are \(C^1\) with
uniformly bounded derivatives.  Differentiating the finite-measure integral
equation with respect to its initial Lamperti coordinate and then using
\(F_\theta'(y)=1/\sigma_\theta(y)\) and
\(G_\theta'(z)=\sigma_\theta(G_\theta(z))\) gives
\eqref{eq:flow-derivative}.  Its exponent is bounded in absolute value by
\(K(T+A)\); the two-sided bounds on \(\sigma_\theta\) then give
\eqref{eq:flow-bi-Lipschitz} and complete the diffeomorphism assertion.
\end{proof}

The following estimate quantifies Helly--Bray continuity on the capacity
cores.

\begin{lemma}[Quantitative Stieltjes stability]
\label{lem:quantitative-stieltjes}
Fix \(0<\eta<1\).  Let \(q,\widetilde q\in\mathcal A_\Lambda\), and let
\(f\in C([0,T])\) satisfy
\[
  \norm{f}_\infty\le B,
  \qquad
  \abs{f(t)-f(s)}\le H\abs{t-s}^\eta.
\]
If \(\delta=\norm{q-\widetilde q}_\infty\), then
\begin{equation}\label{eq:quantitative-stieltjes}
  \sup_{0\le s\le t\le T}
  \left|
    \int_s^t f(r)\,\dd q_r
    -\int_s^t f(r)\,\dd\widetilde q_r
  \right|
  \le C_{B,H,\eta,\Lambda,T}
  \bigl(\delta+\delta^\eta\bigr).
\end{equation}
\end{lemma}

\begin{proof}
The assertion is immediate when \(\delta=0\).  Suppose first that
\(0<\delta\le1\).  For fixed \(s<t\), divide \([s,t]\) into at most
\(1+T/\delta\) intervals of length at most \(\delta\), and let \(f^\delta\)
be the left-step approximation to \(f\) on this mesh.  Since both integrators
have total variation at most \(\Lambda T\),
\[
  \left|
    \int_s^t(f-f^\delta)\,\dd(q-\widetilde q)
  \right|
  \le2\Lambda T H\delta^\eta.
\]
Writing \(h=q-\widetilde q\), discrete summation by parts gives
\[
  \left|\int_s^t f^\delta\,\dd h\right|
  \le2B\delta
     +H\delta\left(1+\frac{T}{\delta}\right)\delta^\eta
  \le C_{B,H,\eta,T}(\delta+\delta^\eta).
\]
These estimates are independent of \(s,t\).  If \(\delta>1\), the direct
bound \(2\Lambda TB\) is absorbed by the right-hand side of
\eqref{eq:quantitative-stieltjes}.
\end{proof}

We use F\"ollmer's deterministic It\^o formula
\cite{Follmer81}.  For \(s\in[0,T]\), let \(\pi_n^s\) be the partition of
\([s,T]\) consisting of \(s\), the dyadic points in \((s,T)\), and \(T\).
Thus every interval \([u,v]\in\pi_n^s\) has \(u\ge s\).  Inserting the single
point \(s\) into the dyadic cell that contains it does not change the
quadratic-variation limit because the oscillation of a continuous path on a
dyadic interval tends to zero.

\begin{lemma}[Uniform weighted quadratic variation on the cores]
\label{lem:uniform-weighted-qv}
Fix \(R,m<\infty\).  Uniformly over
\[
  x\in\cK_R,\quad \theta\in\Theta,\quad
  s\in[0,T],\quad \abs{y}\le m,
\]
the quadratic-variation sums of
\(Z^{\theta,s,y}(x,Q(x))\) along \(\pi_n^s\) define positive discrete
measures on \([s,T]\) which converge weakly to the restriction of
\(\dd Q(x)\) to \([s,T]\), uniformly over the indicated parameters.  More
precisely, if \(\mathscr F\) is any uniformly bounded and equicontinuous
family of continuous weights, then
\begin{equation}\label{eq:uniform-weighted-qv}
  \sup
  \left|
    \sum_{\substack{[u,v]\in\pi_n^s\\v\le t}}
      f(u)\bigl(Z_v-Z_u\bigr)^2
    -
    \int_s^t f(r)\,\dd Q_r(x)
  \right|
  \longrightarrow0,
\end{equation}
where the supremum is also over \(t\in[s,T]\) and \(f\in\mathscr F\).
Furthermore,
\begin{equation}\label{eq:uniform-cubic-remainder}
  \sup
  \sum_{[u,v]\in\pi_n^s}
  \abs{Z_v-Z_u}^3
  \longrightarrow0.
\end{equation}
\end{lemma}

\begin{proof}
The continuous path \(Q^n\) in \eqref{eq:dyadic-qv} is not increasing inside
each grid cell, so we use the positive discrete measure
\[
  \mu_x^n
  =
  \sum_{k=0}^{2^n-1}
  \bigl(x_{t_{k+1}^n}-x_{t_k^n}\bigr)^2
  \delta_{t_{k+1}^n}.
\]
Its cumulative function \(\overline Q^n\) satisfies
\[
  \norm{\overline Q^n(x)-Q^n(x)}_\infty
  \le\omega_x(T2^{-n})^2
  \le R^2\Lambda T2^{-2\eta n},
  \qquad x\in\cK_R.
\]
Together with \eqref{eq:qv-core-rate}, this shows that the cumulative
functions of \(\mu_x^n\) converge to \(Q(x)\), uniformly over
\(x\in\cK_R\).

Restricting the grid to \([s,T]\) and inserting \(s\) changes only the cell
that contains \(s\).  The difference of the two quadratic sums is bounded by
a fixed multiple of the squared oscillation on that cell.  At the moving
upper endpoint there is at most one unfinished cell.  Since
\(Q(x)\in\mathcal A_\Lambda\), the cumulative functions
\[
  M_{n,s}^x(t)
  =
  \sum_{\substack{[u,v]\in\pi_n^s\\v\le t}}
  (x_v-x_u)^2
\]
obey
\begin{equation}\label{eq:inserted-point-qv-rate}
  \sup_{\substack{x\in\cK_R\\0\le s\le t\le T}}
  \abs{
    M_{n,s}^x(t)-\bigl(Q_t(x)-Q_s(x)\bigr)
  }
  \le
  C_{\eta,\gamma}\Lambda T
  \left(
    R^2 2^{-\gamma n}
    +R^2 2^{-2\eta n}
    +2^{-n}
  \right).
\end{equation}

Write \(Z=Z_s+x-x_s+A\), where \(A\) is the finite-variation part in
\eqref{eq:Lamperti-equation}.  On the parameter set in the statement,
put
\[
  B_b=\sup_{\theta\in\Theta}\norm{\widehat b_\theta}_\infty,
  \qquad
  B_c=\sup_{\theta\in\Theta}\norm{\widehat c_\theta}_\infty.
\]
Since \(Q(x)\in\mathcal A_\Lambda\),
\begin{equation}\label{eq:A-variation-modulus}
  \Var(A;[s,T])\le B_bT+B_c\Lambda T,
  \qquad
  \omega_A(\delta)\le(B_b+B_c\Lambda)\delta.
\end{equation}
Hence
\begin{align*}
  \sup_t
  \left|
    \sum_{\substack{[u,v]\in\pi_n^s\\v\le t}}
      \bigl((Z_v-Z_u)^2-(x_v-x_u)^2\bigr)
  \right|
  &\le
  2\omega_x(T2^{-n})\Var(A)
  +\omega_A(T2^{-n})\Var(A)
  \longrightarrow0
\end{align*}
uniformly.  Thus the positive quadratic-variation measures of \(Z\) have
cumulative functions converging uniformly to \(Q_t(x)-Q_s(x)\).

For the endpoint convention, set
\[
  \mu_{n,s}^Z
  =\sum_{[u,v]\in\pi_n^s}(Z_v-Z_u)^2\delta_v.
\]
If \(\omega_{\mathscr F}\) is a common modulus of continuity for
\(\mathscr F\), then
\begin{equation}\label{eq:left-right-weight-error}
  \left|
    \sum_{\substack{[u,v]\in\pi_n^s\\v\le t}}
      f(u)(Z_v-Z_u)^2
    -\int_{[s,t]}f(r)\,\dd\mu_{n,s}^Z(r)
  \right|
  \le
  \omega_{\mathscr F}(T2^{-n})\,
  \mu_{n,s}^Z([s,T]).
\end{equation}
The total masses on the right are uniformly bounded.  Let \(M\) be a common
bound for the masses of \(\mu_{n,s}^Z\) and \(\dd Q(x)|_{[s,T]}\), and let
\(\varepsilon_n\) be the supremum of the cumulative-function errors just
obtained.  Thus \(\varepsilon_n\to0\).  If a fixed deterministic mesh has
mesh size at most \(\rho\) and \(N_\rho\) cells, step-function approximation
and summation by parts give, uniformly in all the parameters and in
\(t\in[s,T]\),
\[
  \left|
    \int_{[s,t]}f(r)\,\dd\mu_{n,s}^Z(r)
    -\int_s^t f(r)\,\dd Q_r(x)
  \right|
  \le
  2M\omega_{\mathscr F}(\rho)
  +2\sup_{g\in\mathscr F}\norm{g}_\infty
      (N_\rho+2)\varepsilon_n.
\]
First let \(n\to\infty\), and then let \(\rho\downarrow0\).  Together with
\eqref{eq:left-right-weight-error}, this proves the weighted convergence in
\eqref{eq:uniform-weighted-qv}.  Finally,
\[
  \sum\abs{\Delta Z}^3
  \le
  \max\abs{\Delta Z}\sum\abs{\Delta Z}^2.
\]
The first factor tends to zero uniformly by equicontinuity, while the second
is uniformly bounded by the already proved quadratic-variation convergence.
This proves \eqref{eq:uniform-cubic-remainder}.
\end{proof}

\begin{theorem}[Simultaneous driver-only It\^o flow]
\label{thm:simultaneous-flow}
For \(x\in\cG_{\mathrm{dyad}}\), put
\begin{equation}\label{eq:driver-flow}
  S_{s,t}^\theta(x;y)
  :=\Phi_{s,t}^\theta(x,Q(x);y).
\end{equation}
When \(S(X)\) is viewed as a random field on all of \(\Omega\), we extend it
by the identity flow \(S_{s,t}^\theta(x;y)=y\) on
\(\cG_{\mathrm{dyad}}^c\).  This yields a total Borel random field.
Then the following statements hold.

\begin{enumerate}[label=(\roman*),leftmargin=2.2em]
  \item For every \(x\in\cG_{\mathrm{dyad}}\), every
  \((\theta,s,y)\in\Theta\times[0,T]\times\R\), and every \(t\ge s\),
  the left sums
  \begin{equation}\label{eq:Follmer-sums}
    J_{s,t}^{n,\theta}(x;y)
    =
    \sum_{[u,v]\in\pi_n^s\,:\,s\le u<t}
    \sigma_\theta\bigl(S_{s,u}^\theta(x;y)\bigr)
    \bigl(x_{v\wedge t}-x_u\bigr)
  \end{equation}
  converge uniformly in \(t\).  Calling the limit
  \(J_{s,t}^\theta(x;y)\), one has the deterministic It\^o equation
  \begin{equation}\label{eq:pathwise-sde-flow}
    S_{s,t}^\theta(x;y)
    =y+
    \int_s^t b_\theta\bigl(S_{s,r}^\theta(x;y)\bigr)\,\dd r
    +J_{s,t}^\theta(x;y).
  \end{equation}
  The convergence and the identity hold for all parameters on the single
  domain \(\cG_{\mathrm{dyad}}\).

  \item The field
  \begin{equation}\label{eq:flow-valued-map}
    x\longmapsto
    \left(
      (\theta,s,t,y)\longmapsto S_{s,t}^\theta(x;y)
    \right)
  \end{equation}
  and the corresponding integral field
  \[
    x\longmapsto
    \left(
      (\theta,s,t,y)\longmapsto J_{s,t}^\theta(x;y)
    \right)
  \]
  are Borel and causal from \(\cG_{\mathrm{dyad}}\) into
  \[
    C_{\mathrm{loc}}
    \bigl(\Theta\times\{(s,t):0\le s\le t\le T\}\times\R;\R\bigr).
  \]
  Here \(C_{\mathrm{loc}}\) carries the compact-open topology.
  The solution field satisfies the cocycle identity \eqref{eq:cocycle} and is a
  \(C^1\)-diffeomorphic flow in the initial value.  The integral field
  satisfies, for \(s\le u\le t\),
  \begin{equation}\label{eq:integral-cocycle}
    J_{s,t}^\theta(x;y)
    =J_{s,u}^\theta(x;y)
     +J_{u,t}^\theta
       \bigl(x;S_{s,u}^\theta(x;y)\bigr).
  \end{equation}

  \item For the compact sets \(\cK_R\) in
  \Cref{cor:parabolic-core}, the flow-valued map
  \eqref{eq:flow-valued-map} and the integral-valued map are uniformly
  continuous on \(\cK_R\).  In fact, with
  \begin{equation}\label{eq:flow-holder-exponent}
    \beta_{\mathrm{flow}}
    =\eta\alpha
    =\frac{\eta\gamma}{1-\eta+\gamma},
  \end{equation}
  there is \(C_R<\infty\) such that, whenever
  \(x,\widetilde x\in\cK_R\) and
  \(e=\norm{x-\widetilde x}_\infty\le1\),
  \begin{align}
    \sup_{\substack{\theta\in\Theta,\ 0\le s\le t\le T\\y\in\R}}
    \Bigl(&
      \abs{S_{s,t}^\theta(x;y)-S_{s,t}^\theta(\widetilde x;y)}
      \notag\\[-0.2em]
      &+
      \abs{J_{s,t}^\theta(x;y)-J_{s,t}^\theta(\widetilde x;y)}
    \Bigr)
    \le C_Re^{\beta_{\mathrm{flow}}}.
    \label{eq:flow-raw-holder}
  \end{align}
  In addition, the sums \eqref{eq:Follmer-sums} converge locally uniformly in
  \((x,\theta,s,t,y)\) on
  \(\cK_R\times\Theta\times\{s\le t\}\times\R\).

  \item Under every \(P\in\mathfrak M_\Lambda\), for every fixed
  \((\theta,s,y)\), the process
  \(S_{s,\cdot}^\theta(X;y)\) is the unique strong solution of
  \begin{equation}\label{eq:classical-scalar-SDE}
    Y_t
    =y+\int_s^t b_\theta(Y_r)\,\dd r
       +\int_s^t\sigma_\theta(Y_r)\,\dd X_r,
    \qquad t\ge s.
  \end{equation}
  For every fixed \((\theta,s,y)\), the corresponding section of \(J\) is a
  version of the It\^o integral in \eqref{eq:classical-scalar-SDE}.  Its
  pathwise equation and cocycle identities hold on
  \(\cG_{\mathrm{dyad}}\) for all parameters, while the stochastic
  identifications hold up to indistinguishability for each fixed parameter.
\end{enumerate}
\end{theorem}

\begin{proof}
Fix \(x\in\cG_{\mathrm{dyad}}\), write \(q=Q(x)\), and abbreviate
\(Z=Z^{\theta,s,y}(x,q)\), \(Y=G_\theta(Z)\).  From
\eqref{eq:Lamperti-equation}, put
\[
  Z_t=Z_s+x_t-x_s+A_t,
\]
where \(A_s=0\) and \(A\) is continuous with finite variation.  Let
\(\Delta_n=T2^{-n}\), and, for a continuous path \(w\), define
\[
  V_{n,s}^w(t)
  =
  \sum_{\substack{[u,v]\in\pi_n^s\\s\le u<t}}
    \bigl(w_{v\wedge t}-w_u\bigr)^2,
  \qquad s\le t\le T.
\]
Comparing the partition obtained by inserting \(s\) with the original
dyadic partition gives
\begin{equation}\label{eq:all-path-inserted-qv}
  \sup_{s\le t\le T}
  \left|
    V_{n,s}^x(t)-(q_t-q_s)
  \right|
  \le
  2\norm{Q^n(x)-q}_\infty+4\omega_x(\Delta_n)^2
  \longrightarrow0.
\end{equation}
Here \(\omega_x\) is the modulus of continuity of \(x\).  Thus
\eqref{eq:all-path-inserted-qv} holds for every
\(x\in\cG_{\mathrm{dyad}}\), independently of the compact cores.

The finite-variation terms satisfy, uniformly in \(t\in[s,T]\),
\[
  \left|V_{n,s}^Z(t)-V_{n,s}^x(t)\right|
  \le
  2\omega_x(\Delta_n)\Var(A;[s,T])
  +\omega_A(\Delta_n)\Var(A;[s,T])
  \longrightarrow0.
\]
Consequently \(V_{n,s}^Z\to q-q_s\) uniformly.  Moreover,
\[
  \sup_{s\le t\le T}
  \left|
    \mu_{n,s}^Z([s,t])-V_{n,s}^Z(t)
  \right|
  \le \omega_Z(\Delta_n)^2\longrightarrow0.
\]
Therefore the positive measures
\[
  \mu_{n,s}^Z
  =\sum_{[u,v]\in\pi_n^s}(Z_v-Z_u)^2\delta_v
\]
have uniformly bounded mass and cumulative functions converging uniformly to
\(q-q_s\).  Step-function approximation of the continuous weight
\(G_\theta''(Z)\), together with the vanishing unfinished-cell error, yields
\begin{equation}\label{eq:all-path-weighted-qv}
  \sup_{s\le t\le T}
  \left|
    \sum_{\substack{[u,v]\in\pi_n^s\\u<t}}
      G_\theta''(Z_u)(Z_{v\wedge t}-Z_u)^2
    -\int_s^tG_\theta''(Z_r)\,\dd q_r
  \right|
  \longrightarrow0.
\end{equation}

Taylor's formula on the intervals of \(\pi_n^s\) gives a remainder
\(R_n(t)\) satisfying, by \eqref{eq:G-third-bound},
\[
  \sup_{s\le t\le T}\abs{R_n(t)}
  \le
  C\omega_Z(\Delta_n)
   \sup_{s\le t\le T}V_{n,s}^Z(t)
  \longrightarrow0.
\]
Together with \eqref{eq:all-path-weighted-qv}, this proves uniform convergence
of the F\"ollmer left sums of \(G_\theta'(Z)\) against \(Z\).  The ordinary
left sums against \(A\) also converge uniformly, since \(A\) has finite
variation.  Subtracting them leaves precisely the sums
\eqref{eq:Follmer-sums}.  Their limit is
\begin{align*}
  Y_t-y
  &-\int_s^tG_\theta'(Z_r)\widehat b_\theta(Z_r)\,\dd r\\
  &-\int_s^t
    \left(
      G_\theta'(Z_r)\widehat c_\theta(Z_r)
      +\frac12G_\theta''(Z_r)
    \right)\dd q_r.
\end{align*}
Since
\[
  G_\theta'(Z_r)\widehat b_\theta(Z_r)
    =b_\theta(Y_r),
  \qquad
  G_\theta'(Z_r)\widehat c_\theta(Z_r)
    +\frac12G_\theta''(Z_r)=0,
\]
the limit equals \(Y_t-y-\int_s^t b_\theta(Y_r)\,\dd r\).  This proves the
uniform-in-\(t\) convergence in (i) and \eqref{eq:pathwise-sde-flow} for all
parameters on \(\cG_{\mathrm{dyad}}\).

Set
\[
  \Delta_T=\{(s,t):0\le s\le t\le T\},
  \qquad
  \mathsf P_0=\Theta\times\Delta_T\times\R,
\]
and exhaust \(\mathsf P_0\) by the compact sets
\(\mathsf P_m=\Theta\times\Delta_T\times[-m,m]\).  On every set
\(\{q_T\le a\}\), \Cref{lem:deterministic-flow} and compactness of
\(\mathsf P_m\) show that
\[
  (x,q)\longmapsto
  \left(
    (\theta,s,t,y)\longmapsto\Phi_{s,t}^\theta(x,q;y)
  \right)
\]
is continuous into \(C_{\mathrm{loc}}(\mathsf P_0;\R)\).  Since
\[
  C_\uparrow([0,T])
  =\bigcup_{a\in\N}\{q\in C_\uparrow([0,T]):q_T\le a\}
\]
is a countable union of closed sets, the same function-space-valued map is
Borel on \(\Omega\times C_\uparrow([0,T])\).  Composing it with the Borel
map \(x\mapsto(x,Q(x))\) proves that the solution field in
\eqref{eq:flow-valued-map} is a
\(C_{\mathrm{loc}}(\mathsf P_0;\R)\)-valued Borel map.  Causality follows
from the causality of \(Q\) and of the deterministic equation.
Because \(\cG_{\mathrm{dyad}}\) is Borel and the identity flow is a fixed
element of the target space, the extension on \(\cG_{\mathrm{dyad}}^c\) is
Borel as well.

Equation \eqref{eq:pathwise-sde-flow} gives the identity
\begin{equation}\label{eq:J-from-S}
  J_{s,t}^\theta(x;y)
  =S_{s,t}^\theta(x;y)-y
   -\int_s^t b_\theta(S_{s,r}^\theta(x;y))\,\dd r.
\end{equation}
The operator on \(C_{\mathrm{loc}}(\mathsf P_0;\R)\) defined by the
right-hand side of \eqref{eq:J-from-S} is continuous; on each
\(\mathsf P_m\), its Lipschitz constant is at most \(1+LT\).  Hence the
integral field is also Borel and causal.  Splitting
\eqref{eq:J-from-S} at \(u\) and using the solution cocycle proves
\eqref{eq:integral-cocycle}.  Formula
\eqref{eq:flow-derivative} gives the \(C^1\)-diffeomorphism assertion and
completes (ii).

On \(\cK_R\), the map
\(x\mapsto(x,Q(x))\) is continuous and has compact image.  Composition with
\eqref{eq:joint-flow-continuity} is therefore continuous into the stated
local-uniform function space and hence uniformly continuous; the same is
true for \(J\) by \eqref{eq:J-from-S}.

We prove the stronger modulus \eqref{eq:flow-raw-holder}.  Take
\(x,\widetilde x\in\cK_R\), put
\(q=Q(x)\), \(\widetilde q=Q(\widetilde x)\), and write
\(\delta=\norm{q-\widetilde q}_\infty\).  For fixed
\((\theta,s,y)\), let \(Z,\widetilde Z\) denote the two Lamperti solutions.
Write again
\[
  B_b=\sup_{\theta\in\Theta}\norm{\widehat b_\theta}_\infty,
  \qquad
  B_c=\sup_{\theta\in\Theta}\norm{\widehat c_\theta}_\infty.
\]
The boundedness of \(\widehat b,\widehat c\), the \(\eta\)-H\"older bound on
\(\widetilde x\), and \(\widetilde q\in\mathcal A_\Lambda\) give, uniformly
in all parameters,
\[
  \abs{\widetilde Z_t-\widetilde Z_u}
  \le R\sqrt\Lambda\,T^{1/2-\eta}\abs{t-u}^\eta
      +(B_b+B_c\Lambda)\abs{t-u}
  \le C_R\abs{t-u}^\eta.
\]
Thus \(r\mapsto\widehat c_\theta(\widetilde Z_r)\) belongs to a uniformly
bounded, uniformly \(\eta\)-H\"older family.  Subtracting the two Lamperti
equations, applying \Cref{lem:quantitative-stieltjes} to the term driven by
\(q-\widetilde q\), and then applying Gronwall's inequality for the measure
\(\dd r+\dd q_r\), yields
\[
  \sup_{t\in[s,T]}\abs{Z_t-\widetilde Z_t}
  \le C_R\bigl(e+\delta+\delta^\eta\bigr),
  \qquad e=\norm{x-\widetilde x}_\infty.
\]
Since the maps \(G_\theta\) are uniformly Lipschitz and
\(\delta\le C_Re^\alpha\) by \eqref{eq:qv-core-modulus}, we obtain
\[
  \sup_{\substack{\theta\in\Theta,\ 0\le s\le t\le T\\y\in\R}}
  \abs{S_{s,t}^\theta(x;y)-S_{s,t}^\theta(\widetilde x;y)}
  \le C_Re^{\eta\alpha},
  \qquad e\le1.
\]
The same estimate for \(J\) follows from \eqref{eq:J-from-S} and the common
Lipschitz bound on \(b_\theta\).  This proves \eqref{eq:flow-raw-holder}.

The local uniform assertion for the sums follows
from the quantitative deterministic proof of F\"ollmer's formula.  On
\(\cK_R\times\Theta\times\{s\le t\}\times[-m,m]\), the families \(x\),
\(q=Q(x)\), \(Z\), and their finite-variation parts are uniformly bounded
and equicontinuous; the latter also have uniformly bounded variation.
Taylor's formula separates four errors.  First, the finite-variation left
sums converge uniformly because
\[
  \sup_t\left|
    \sum_{\substack{[u,v]\in\pi_n^s\\u<t}}
      G_\theta'(Z_u)(A_{v\wedge t}-A_u)
    -\int_s^tG_\theta'(Z_r)\,\dd A_r
  \right|
  \le
  \omega_{G_\theta'(Z)}(T2^{-n})\Var(A;[s,T])
  \longrightarrow0.
\]
Second, the uniform third-derivative bound
\eqref{eq:G-third-bound} and
\Cref{lem:uniform-weighted-qv}, applied to the equicontinuous family of
weights \(G_\theta''(Z)\), give uniform convergence of the quadratic terms.
Third, the cubic Taylor remainders vanish by
\eqref{eq:uniform-cubic-remainder}.  Fourth, the sole unfinished terminal
cell is controlled by the common oscillation modulus of \(Z\).  These four
estimates prove the local uniform convergence in (iii).

Finally fix \(P\in\mathfrak M_\Lambda\) and an adapted continuous version of
\([X]^P\).  Solve \eqref{eq:Lamperti-equation} with the two adapted inputs
\((X,[X]^P)\).  Causality of the deterministic equation, or equivalently its
Picard iteration, shows that, for every fixed \((\theta,s,y)\), the resulting
processes \(\overline Z\) and \(\overline S=G_\theta(\overline Z)\) are
adapted.  By
\eqref{eq:intrinsic-qv-identification},
\[
  \overline S_{s,\cdot}^\theta(y)
  =S_{s,\cdot}^\theta(X;y)
  \qquad P\text{-almost surely}.
\]
The equality holds simultaneously in the parameters on the single
\(P\)-full event \(\{Q(X)=[X]^P\}\).  Under the usual augmentation, its
complement belongs to \(\cF_0^P\).  Hence, for every fixed \((\theta,s,y)\),
the section \(S_{s,\cdot}^\theta(X;y)\) is adapted under \(P\); continuity
then implies progressive measurability.  The process
\(\overline Z\) is a continuous semimartingale.  Applying the classical
It\^o formula to \(G_\theta(\overline Z)\) and using the two cancellation
identities proves \eqref{eq:classical-scalar-SDE}.  Global Lipschitz
continuity gives strong uniqueness.  Equation \eqref{eq:J-from-S} then
identifies \(J\), up to indistinguishability, with the It\^o integral for each
fixed \((\theta,s,y)\), proving (iv).
\end{proof}

\begin{corollary}[A simultaneous stochastic flow]
\label{cor:simultaneous-stochastic-flow}
The fixed Borel set
\(N=\Omega\setminus\cG_{\mathrm{dyad}}\) is
\(\mathfrak M_\Lambda\)-polar.  On \(N^c\), simultaneously for every
\(\theta\in\Theta\), every starting time \(s\), and every initial value
\(y\), the fields \((S,J)\) satisfy the pathwise It\^o equation, the solution
and integral cocycles, and form an orientation-preserving \(C^1\) flow in
\(y\).  Under each \(P\in\mathfrak M_\Lambda\), every fixed parameter section is the
unique classical strong solution and its It\^o integral version.
\end{corollary}

\begin{proof}
Combine \eqref{eq:Gdyad-full-capacity} with
\Cref{thm:simultaneous-flow}.
\end{proof}

\begin{remark}[Comparison with endogenous-partition sums]
For every fixed parameter, the Bichteler--Karandikar integral from
\Cref{sec:integral} and the F\"ollmer integral in
\Cref{thm:simultaneous-flow} agree quasi surely.  Uniform stability of the
endogenous hitting-time sums over the full parameter class remains open
because their crossing times depend recursively on the integrand.
\end{remark}

\section{Pathwise geometry and capacity cores}
\label{sec:geometry}

\subsection{Relation to rough-path integration}

The endogenous functional \(\cI\) is built from the order of time and the
level crossings of the integrand.  Its stopping rule
\eqref{eq:tau-recursion}, sums \eqref{eq:In-def}, and convergence domain
\eqref{eq:Dt-def} are deterministic.  The semimartingale identification
theorem supplies convergence and identifies the pathwise limit under each
law.

Rough-path integration begins with an enhanced driver, such as
\(\mathbf x=(x,\mathbb X)\), and obtains continuity in a rough-path topology.
The endogenous construction instead uses a partition generated by \(h\), so
its natural input is the pair \((h,x)\); its structural properties are Borel
measurability, projective consistency, and causality on a partial domain.
Under semimartingale laws it produces the It\^o integral directly.  In a
rough-path formulation the same integral is represented by a non-geometric
It\^o enhancement, or by a geometric enhancement with the bracket
correction.  Whenever both constructions represent the same stochastic
integral on a common canonical space, \Cref{cor:scheme-independence} gives
their quasi-sure equality.

The driver-only construction in \Cref{thm:simultaneous-flow} uses a different
mechanism: continuous dyadic quadratic variation supplies the second-order
input to the Lamperti--F\"ollmer equation.  Its continuity is measured on the
raw-path cores \(\cK_{A,B}\), rather than in a prescribed rough-path
topology.

\subsection{Stable capacity cores}

A common intrinsic domain gives simultaneous existence.  Quantitative
stability is encoded by subsets on which the solution map has a deterministic
modulus of continuity.

\begin{definition}[Stable capacity core]\label{def:capacity-core}
Let \(S:G\to F\) be a solution map on a Borel subset of a Polish input space
\((E,d_E)\), with \(F\) metrized by \(d_F\).  A sequence of Borel sets
\((K_R)_{R\ge1}\) is a stable capacity core if
\[
  K_R\subset K_{R+1}\subset G,
  \qquad
  c_{\cP}(K_R^c)\longrightarrow0,
\]
and, for every \(R\), there is a modulus
\(\omega_R:[0,\infty)\to[0,\infty)\), with
\(\omega_R(r)\downarrow0\) as \(r\downarrow0\), such that
\[
  d_F(S(x),S(y))
  \le
  \omega_R(d_E(x,y)),
  \qquad x,y\in K_R.
\]
\end{definition}

On stable capacity cores, approximation and parameter-uniform estimates
reduce to deterministic continuity bounds.

\Cref{thm:driver-core,cor:parabolic-core,thm:simultaneous-flow} provide a
concrete solution to this problem for scalar, uniformly positive autonomous
equations
\eqref{eq:classical-scalar-SDE}, uniformly over the compact coefficient
family \(\Theta\).  They also give stopping and concatenation stability of
the driver cores.  This leads to the following extension problem.

\begin{problem}[General capacity-core solution principle]
\label{problem:capacity-core}
Identify verifiable conditions on a nondominated semimartingale family
\(\cP\), a coefficient class, and a canonical approximation scheme such that
there exist:
\begin{enumerate}[label=(\roman*),leftmargin=2.2em]
  \item a Borel causal solution map \(S\) for the corresponding law-free
  relation;
  \item stable capacity cores \(K_R\) with
  \(c_{\cP}(K_R^c)\to0\);
  \item estimates on \(K_R\) uniform in the initial value, starting time, and
  coefficient parameter.
\end{enumerate}
The desired conclusion is a single pathwise flow, with the flow and cocycle
identities holding outside one polar set simultaneously for the entire
parameter class.
\end{problem}

\begin{remark}
Li and Scheutzow \cite{LiScheutzow11} give smooth bounded-coefficient
examples in which each fixed initial condition has a global solution almost
surely, while a globally defined continuous flow for all initial conditions
fails.  Uniform pathwise estimates are therefore an additional requirement
for the general capacity-core problem.
\end{remark}

The input space also depends on the model class.  For the all-Dirac family in
\cite[Example 26]{DenisHuPeng11}, the only polar set is empty, so a common
domain must contain every path.  Singular equations may require an enhanced
or renormalized input; that structure then becomes part of the canonical
state space.

\subsection{Conclusion}

The endogenous-partition sums define a Borel causal It\^o functional on one
fixed pair-domain of full measure under every continuous-semimartingale law.
For the maximal bounded-volatility martingale class, the dyadic construction
produces two-scale compact driver cores with explicit capacity tails,
approximation rates, and raw-path continuity of the intrinsic quadratic
variation.

The Lamperti--F\"ollmer equation lifts these estimates to one Borel causal
solution field, simultaneously in the coefficient, starting time, and initial
value.  Outside a single polar set, this field is an orientation-preserving
\(C^1\) flow and satisfies the solution and integral cocycles.  Extensions to
multidimensional, degenerate, and singular equations require corresponding
driver-only cores and canonical second-order inputs.

\bibliographystyle{alpha}

\end{document}